\documentclass[12pt]{amsart}
\usepackage{amssymb}
\usepackage{graphicx}
\input epsf

\newcommand {\supplus}{\mathop{{\supset}\llap{\raise
0.5pt\hbox{\normalfont\small+}\hskip 0.5pt}}}

\newcommand {\subplus}{\mathop{{\subset}\llap{\raise
0.5pt\hbox{\normalfont\small+}\hskip 0.5pt}}}

\newcommand {\Cee}    {{\mathbb  C}}

\newcommand {\Nee}    {{\mathbb  N}}

\newcommand {\Zee}    {{\mathbb  Z}}

\newcommand {\fa}     {{\mathfrak{a}}}
\newcommand {\fA}     {{\mathfrak{A}}}

\newcommand {\faut}   {{\mathfrak{aut}}}
\newcommand {\fb}     {{\mathfrak{b}}}

\newcommand {\fg}     {{\mathfrak{g}}}    %
\newcommand {\fgl}    {{\mathfrak{gl}}}  %
\newcommand {\fh}     {{\mathfrak{h}}}

\newcommand {\fL}     {{\mathfrak{L}}}

\newcommand {\fM}     {{\mathfrak{M}}}   %

\newcommand {\fo}     {{\mathfrak{o}}}
\newcommand {\fosp}   {{\mathfrak{osp}}}
\newcommand {\fpe}    {{\mathfrak{pe}}}   %
\newcommand {\fpg}    {{\mathfrak{pg}}}

\newcommand {\fpgl}   {{\mathfrak{pgl}}}

\newcommand {\fpsl}   {{\mathfrak{psl}}}
\newcommand {\fpq}    {{\mathfrak{pq}}}

\newcommand {\fpsq}   {{\mathfrak{psq}}}
\newcommand {\fq}     {{\mathfrak{q}}}

\newcommand {\fs}     {{\mathfrak{s}}}
\newcommand {\fS}     {{\mathfrak{S}}}

\newcommand {\fsh}    {{\mathfrak{sh}}}

\newcommand {\fsl}    {{\mathfrak{sl}}}

\newcommand {\fsp}    {{\mathfrak{sp}}}
\newcommand {\fspe}   {{\mathfrak{spe}}}

\newcommand {\fsq}    {{\mathfrak{sq}}}

\newcommand {\cal} {\mathcal}

\newcommand {\cA}     {{\cal A}}

\newcommand {\cL}     {{\cal L}}

\newcommand {\cO}     {{\cal O}}

\def \opname#1#2%
  {\expandafter\newcommand \csname #1\endcsname {{\mathop{#2}\nolimits}}}

\newcommand{\rmname}[1]
  {\expandafter\newcommand \csname #1\endcsname {{\operatorname{#1}}}}

\newcommand{\rmnameii}[2]
  {\expandafter\newcommand \csname #1\endcsname {{\operatorname{#2}}}}

\rmname{act}
\rmname{Ad}
\rmname{Add}
\rmname{ad}
\rmname{Alt}
\rmname{alt}
\rmname{Ann}
\rmname{antidiag}
\rmname{Ber}
\rmname{ber}
\rmname{Br}
\rmname{card}
\rmname{ch}
\rmname{Char}
\rmname{cem}
\rmname{cj}
\rmname{Cliff}
\rmname{cntr}
\rmname{codim}
\rmname{coind}
\rmname{const}
\rmname{col}
\rmname{cork}
\rmname{cpr}
\rmname{diag}
\rmnameii{Div}{div}
\rmname{Def}
\rmname{Der}
\rmname{Dim}
\rmname{End}
\rmname{Even}
\rmname{Ext}
\rmname{gr}
\rmname{Hom}
\rmname{HT}
\rmnameii{Ht}{ht}
\rmname{hwt}
\rmname{Id}
\rmname{id}
\rmname{ind}
\rmname{Ind}
\rmname{Inf}
\rmname{irr}
\rmname{Le}
\rmname{Lie}
\rmname{lwt}
\rmname{mult}
\rmname{Mor}
\rmname{nm}
\rmname{Ob}
\rmname{Odd}
\rmname{Osc}
\rmname{per}
\rmname{Pic}
\rmname{pr}
\rmname{pro}
\rmname{Prime}
\rmname{Proj}
\rmname{prt}
\rmname{pt}
\rmname{Q}
\rmname{qet}
\rmname{qtr}
\rmname{rd}
\rmname{rk}
\rmname{row}
\rmname{Res}
\rmname{salt}
\rmname{Sch}
\rmname{sch}
\rmname{SBr}
\rmname{scalar}
\rmname{Ser}
\rmname{sign}
\rmname{Smbl}
\rmname{spin}
\rmname{ssym}
\rmname{str}
\rmname{st}
\rmname{sgn}
\rmname{sq}
\rmname{symm}
\rmname{supp}
\rmname{Supp}
\rmname{St}
\rmname{Spec}
\rmname{Spm}
\rmname{tr}
\rmname{vpt}
\rmname{weyl}
\rmname{Weyl}
\rmname{Witt}

\opname{vvol}  {{v\hspace{-0.1ex}o\hspace{-0.02ex}l\/}}
\opname{pnt}  {\text{\normalfont pt}}
\opname{Span} {{Span}}
\opname{slim} {\overline{\lim}}
\opname{Vol}  {{V\hspace{-0.55ex}o\hspace{-0.02ex}l\/}}
\opname{Par}  {{P\hspace{-0.3ex}a\hspace{-0.05ex}r\/}}

\rmname{Mat}
\rmname{Bil}
\rmname{Diff}
\rmname{Ker}
\rmname{Herm}
\rmname{Coker}
\rmname{Conn}
\rmname{Covect}
\rmname{Vect}
\rmname{Int}

\rmnameii {IM} {Im}
\rmnameii {RE} {Re}

\opname{Aut} {{A\hspace{-0.2ex}u\hspace{-0.1ex}t\/}}
\opname{GL} {{G\hspace{-0.3ex}L}}
\opname{SL} {{S\hspace{-0.3ex}L}}
\opname{Exp} {{E\hspace{-0.2ex}x\hspace{-0.1ex}p\/}}
\opname{GQ} {{G\hspace{-0.2ex}Q}}
\opname{OSp} {{O\hspace{-0.25ex}S\hspace{-0.15ex}p\/}}
\opname{Out} {{O\hspace{-0.25ex}u\hspace{-0.15ex}t\/}}
\opname{Spp} {{S\hspace{-0.2ex}p\/}}
\opname{SpO} {{S\hspace{-0.2ex}p\hspace{-0.02ex}O\/}}
\opname{Pe} {{P\hspace{-0.25ex}e\/}}
\opname{SPe} {{S\hspace{-0.25ex}P\hspace{-0.25ex}e\/}}
\opname{Spin} {{S\hspace{-0.25ex}p\hspace{-0.05ex}i\hspace{-0.1ex}n\/}}
\opname{Iso} {{I\hspace{-0.25ex}s\hspace{-0.1ex}o\/}}
\opname{SSPe} {{S\hspace{-0.25ex}S\hspace{-0.15ex}P\hspace{-0.25ex}e\/}}
\opname{PeU} {{P\hspace{-0.25ex}e\hspace{-0.1ex}U\/}}
\opname{QU} {{Q\hspace{-0.15ex}U\/}}
\opname{U} {{U\/}}

\opname{cGQ} {{\cal G \hspace{-0.2em} Q \/}}
\opname{cSL} {{\cal S \hspace{-0.2em} L \/}}
\opname{cGL} {{\cal G \hspace{-0.2em} L \/}}
\opname{cOSp} {{\cal O \hspace{-0.2em} S \hspace{-0.3em} \it p\/}}
\opname{cPe} {{\cal P \hspace{-1.5pt} \it e\/}}
\opname{cVect} {{\cal V \hspace{-1.5pt} \it
e\hspace{-0.1ex}c\hspace{-0.1ex}t\/}}
\opname{cVol} {{\cal V \hspace{-1.5pt} \it o\hspace{-0.1ex}l\/}}
\opname{cAut} {{\cal A \hspace{-0.2em} \it u\hspace{-0.1em}t\/}}
\opname{cCovect} {{\cal C \hspace{-1.5pt}
     \it o\hspace{-0.1ex}v\hspace{-0.1ex}e\hspace{-0.1ex}c\hspace{-0.1ex}t\/}}
\opname{CW} {{C\hspace{-0.15ex}W}}

\newcommand {\ev} {{\bar0}}
\newcommand {\od} {{\bar1}}

\newcommand {\tto} {\longrightarrow}

\newcommand {\bcdot}   {\mathbin{\hbox{\raise.4ex\hbox{\bf.}}}} 

\newcommand {\secno} {}
\newcommand {\ssecfont} {\normalfont\bf}

\newtheorem{Theorem}{\secno Theorem}

\newtheorem{Lemma}[Theorem]{\secno Lemma}

\newtheorem{Corollary}[Theorem]{\secno Corollary}

\newenvironment {th*}[1]
    {\gdef\thname{#1} \begin{thn}}%
    {\end{thn}}
\newtheorem{thn}[Theorem] {\thname}

\theoremstyle{definition}

\newenvironment {ex*}[1]
    {\gdef\thname{#1} \begin{exn}}%
    {\end{exn}}
\newtheorem{exn}[Theorem]{\thname}

\theoremstyle{remark}

\newtheorem{Remark}[Theorem]{\secno Remark}

\newenvironment {rem*}[1]
    {\gdef\thname{#1} \begin{remn}}%
    {\end{remn}}
\newtheorem{remn}[Theorem]{\thname}

\newcommand {\ssec}{\subsection*}

\newcommand {\ssbegin}[2]
  {\def \secno {\gdef \secno {}{\ssecfont #1. }}%
   \begin{#2}}
\begin{document}
\title[Classical invariant theory for Lie superlagebras]{An analog of
the classical invariant theory \\ for Lie superlagebras}

\author{Alexander Sergeev}

\address{Dept.  of Math., Univ.  of Stockholm, Roslagsv.  101,
Kr\"aftriket hus 6, S-106 91, Stockholm, Sweden (On leave of absence
from Balakovo Inst.  of Technology Technique and Control)\\ e-mail:
mleites@matematik.su.se subject: for Sergeev}

\thanks{I am thankful to D. Leites for support and help.}

\begin{abstract} Let $V$ be a finite-dimensional superspace over
$\Cee$ and $\fg$ a simple (or a ``close'' to simple) matrix Lie
superalgebra, i.e., a Lie subsuperalgebra in $\fgl(V)$.  Under the
{\it classical invariant theory} for $\fg$ we mean the description of
$\fg$-invariant elements of the algebra
$$
\fA^{p, q}_{k, l}=S^{\bcdot}(V^{k}\oplus \Pi (V)^{l}\oplus V^{*p}\oplus
\Pi (V)^{*q}).
$$
We give such description for $\fgl(V)$, $\fsl(V)$ and $\fosp(V)$ and
their ``odd'' analogs: $\fq(V)$, $\fs\fq(V)$; $\fpe(V)$ and
$\fspe(V)$.
\end{abstract}

\subjclass{17A70, 13A50}

\keywords{Invariant theory, Lie superalgebras.}

\maketitle

This is a detailed exposition of my short announcement ``An Analog
of the Classical Invariant theory for Lie Superalgebras'' published
in Funktsional'nyj Analiz i ego Prilozheniya, 26, no. 3, 1992,
88--90. For prerequisites on Lie superalgebras see Appendix borrowed
from \cite{L}.

\section*{\S 1. Setting of the problem. Formulation of the results}

\ssec{1.0} Let $V$ be a finite-dimensional superspace over $\Cee$ and
$\fg$ an arbitrary matrix Lie superalgebra, i.e., a Lie
subsuperalgebra in $\fgl(V)$. Under the {\it classical invariant
theory} for $\fg$ we mean the description of $\fg$-invariant
elements of the algebra
$$
\fA^{p, q}_{k, l}=S^{\bcdot}(V^{k}\oplus \Pi (V)^{l}\oplus V^{*p}\oplus
\Pi (V)^{*q}).
$$
Clearly, $\fA^{p, q}_{k, l}=S^{\bcdot}(U\otimes V\bigoplus V^{*}\otimes
W)$, where
$\dim U=(k, l)$ and $\dim W=(p, q)$.  Therefore, on $\fA^{p, q}_{k,
l}$ there also act Lie superalgebras $\fgl (U)$ and $\fgl (W)$; hence,
the universal enveloping algebra $U(\fgl (U\otimes W))$ also acts on
$\fA^{p, q}_{k, l}$.  The elements of the enveloping algebra will be
called {\it polarization operators}.  These operators commute with the
natural $\fgl(V)$-action.

A set $\fM$ of $\fg$-invariants will be called a {\it basic} one if
the algebra of invariants coincides with the least subalgebra
containing $\fM$ and invariant with respect to polarization operators.

For every series of the classical Lie superalgebras (i.e., simple ones
and their central extensions) we describe such a set $\fM$.

Introduce $\Zee_{2}$-graded sets:
\begin{align*}
T=\{1, \ldots , k, \od, \dots , \bar e\}, \cr
S=\{1, \dots , p, \od, \dots , \bar q\}, \cr
I=\{1, \dots , n, \od, \dots , \overline m\} \cr
\end{align*}
here the odd elements are bared and even elements not bared and select
bases in the spaces $U$, $W$, $V$:
$$
 \{u_{t}\}_{t\in T}, \quad \{w_{s}\}_{s\in S}, \quad \{l_{i}\}_{i\in I}
$$
so that the parity of the vector of a basis coincides with the parity
of the corresponding index.  In $V^{*}$, select the basis
$\{e^{*}_{i}\}_{i\in I}$ left dual to $\{e_{i}\}_{i\in I}$.  Denote:
$$
x_{ti}=u_{t}\otimes e_{i};\quad x^{*}_{is}=e^{*}_{i}\otimes w_{s},
$$
where $v_{s}$ is the column-vector with coordinates $x^{*}_{1s}, \dots
, x^{*}_{\overline ms}$ and $v^{*}_{t}$ is the row-vector with coordinates
$x_{t1}, \dots , x_{t\overline m}$.  Define the scalar product by setting
$$
(v^{*}_{t}, v_{s})=\sum^{}_{i\in I}x_{ti}x^{*}_{is}.
$$
Set
\begin{align*}
\Delta &=\det (x_{ti})_{t, i\in I_\ev};\quad \Delta ^{*}=
\det (x^{*}_{is})_{i, s\in I_\ev}; \cr
\omega &=\det (x_{ti})_{t, i\in I_{\od}};\quad \omega ^{*}=\det
(x^{*}_{is})_{i, s\in I_{\od}}. \end{align*}

\ssbegin{1.1}{Theorem} For the invariants of $\fgl(V)$ a basic set is
the collection of scalar products
$$
(v^{*}_{t}, v_{s})\; \text{ for }\; t\in T, \; s\in S.
$$
\end{Theorem}

\ssbegin{1.1.1}{Corollary} The scalar products $(v^{*}_{s}, v_{t})$
for $t\in T$, $s\in S$ constitute a system of generators of
$\fgl(V)$-invariants.
\end{Corollary}

\ssbegin{1.2}{Theorem} For the invariants of $\fsl (V)$ a basic set
consists of

{\em a)} basic invariants $\fgl(V)$;

{\em b)} the collection of polynomials for $ k\in \Nee$
\begin{align*}
f_{k}&=(\Delta ^{*})^{k}\omega ^{k}\prod \limits_
{t\in I_{\od}, \; s\in I_\ev}(v^{*}_{t}, v_{s}), \cr
f_{-k}&=\Delta ^{k}(\omega ^{*})^{k}\prod \limits_
{t\in I_\ev, \; s\in I_{\od}}(v^{*}_{t}, v_{s}).
\end{align*}
\end{Theorem}

\ssec{1.3} Let $\fosp (V)$ be the Lie superalgebra preserving
the tensor
$$
\sum^{}_{i\in I_\ev}e^{*}_{i}\otimes e^{*}_{n-i+1}+
\sum^{r}_{j=1}(e^{*}_{\overline{m-j+1}}\otimes
e^{*}_{\bar j}-e^{*}_{\bar j}\otimes e^{*}_{\overline{m-j+1}})\;
\text{ for }\; m=2r.
$$
Then the inner products
$$
(v_{s}, v_{t})=\sum^{}_{i\in I_\ev}x^{*}_{is}x^{*}_{n-i+1t}+
(-1)^{p(s)}\sum^{r}_{j=1}(x^{*}_{\overline{m-j+1}s}
x^{*}_{\bar jt}-x^{*}_{\bar js}x^{*}_{\overline{m-j+1}t})
$$
are $\fosp (V)$-invariant.  In what follows we will show that there
also exists an invariant polynomial $\Omega$ such that
$$
\Omega ^{2}=(\det (v_{s}, v_{t})_{s, t\in I_\ev})^{2r+1}.
$$
The existence of an even $\fosp (V)$-invariant form determines an
isomorphism of algebras and $\fosp (V)$-modules $\fA^{p, q}_{k,
l}=\fA^{p+k, q+l}$.  Therefore, we may (and will) confine ourselves to
the case $k=l=0$.

\begin{Theorem} For the invariants of $\fosp (V)$ a basic set
consists of

{\em a)} the scalar products $(v_{s}, v_{t})$ for $s, t\in S$;

{\em b)} the polynomial $\Omega$.
\end{Theorem}

\ssec{1.4} Let $\dim V=(n, n))$ and $\fpe (V)$ be the Lie superalgebra
preserving the tensor
$$
\sum^{}_{i\in I_\ev}(e^{*}_{i}\otimes
e^{*}_{\bar i}+e^{*}_{\bar i}\otimes e^{*}_{i})
$$
Then the inner products
$$
(v_{s}, v_{t})=\sum^{}_{i\in I_\ev}((-1)^{p(s)}
x^{*}_{is}x^{*}_{\bar it}+x^{*}_{\bar is}x^{*}_{it})
$$
are $\fpe(V)$-invariants.

The existence of an odd $\fpe(V)$-invariant form
determines an isomorphism of algebras and $\fpe(V)$-modules
$\fA^{p, q}_{k, l}=\fA^{p+l, q+k}$. Therefore,
we may (and will) assume that $k=l=0$.

\ssbegin{1.4.1}{Theorem} For the invarians of $\fpe(V)$ a basic
set is constituted by the inner products $(v_{s}, v_{t})$ for $ s, t\in S$.
\end{Theorem}

\ssbegin{1.4.1.1}{Corollary} The inner products form a system of
generators of the algebra of $\fpe(V)$-invariants.
\end{Corollary}

\ssbegin{1.4.2}{Theorem} For the invariants of $\fspe(V)$ a basic set
is formed by

{\em a)} the basic invariants for $\fpe(V)$;

{\em b)} the polynomials
$$
p_{k}=\Delta ^{*k}\prod _{s\le t,\; s, t\in I_\ev}(v_{s}, v_{t})
\quad\text{and}\quad p_{-k}=\omega ^{*k}\prod _{s<t,\; s, t\in
I_{\bar 1}}(v_{s}, v_{t}).
$$
\end{Theorem}

\ssec{1.5} Let $\dim V=(n, n)$ and $\fq (V)$ be the Lie
superalgebra preserving the tensor
$$
\sum^{}_{i\in I_\ev}(e_{i}\otimes e^{*}_{\bar i}+
e_{\bar i}\otimes e^{*}_{i}).
$$
Then the expression
$$
[v^{*}_{t}, v_{s}]=\sum (x_{it}x^{*}_{\bar is}+x_{\bar is}
x^{*}_{is})\; \text {for any} \; t\in T_\ev, \; s\in S_\ev
$$
is a $\fq (V)$-invariant.  Since there is an isomorphism of algebras
$\fA^{p, q}_{k, l}=\fA^{p+q}_{k+l}$ and $\fq (V)$-modules, we may (and
will) assume that $q=l=0$.

\ssbegin{1.5.1}{Theorem} For the invariants of $\fq (V)$ a basic
set is the collection of inner products
$$
(v^{*}_{t}, v_{s}), \; [v^{*}_{t}, v_{s}], \quad \; \text {for
any} \; t\in T_\ev, \; s\in S_\ev
$$
\end{Theorem}

\begin{Corollary} The inner products form a system of
generators of the algebra of $\fq (V)$-invariants.
\end{Corollary}

Let $\fsq (V)$ be the queertraceless subalgebra of $\fq (V)$.
Let $Z$ be a matrix of the form
$$
Z=\begin{pmatrix} Z_{0} & Z_{1} \cr Z_{1} & Z_{0} \end{pmatrix}
$$
where
$$
Z_{0}=\{(v^{*}_{t}, v_{s})\}_{t, s\in I_\ev},
\quad Z_{1}=\{[v^{*}_{t}, v_{s}]\}_{t, s\in I_\ev}
$$
and $Y$ be a matrix of the form
$$
Y=\begin{pmatrix} Y_{0} & Y_{1} \cr Y_{1} & Y_{0} \end{pmatrix},
$$
where
$$
Y_{0}=\{x^{*}_{is}\}_{i, s\in I_\ev}, \quad
Y_{1}=\{x^{*}_{it}\}_{i \in I_\ev, t\in I_{\od}}.
$$
In what follows we will prove that for any partition of $\lambda $
$$
\lambda=\lambda _{1}>\lambda _{2}>\dots
>\lambda _{n}>0
\eqno{(1.5.1)}
$$
the expression
$$
q_{\lambda }={\rm qtr}\; Z^{\lambda _{1}}\cdot \dots\cdot
{\rm qtr}\; Z^{\lambda _{n}}\cdot \qet\; Y
$$
is a polynomial.

\ssbegin{1.5.2}{Theorem} A basic
set for the invariants of $\fsq (V)$ is formed by

{\em a)} the basic invariants for $\fq (V)$ and

{\em b)} the polynomials $q_{\lambda }$, where $\lambda $ runs over
all partitions of the form $(1.5.1)$.
\end{Theorem}

\section*{\S 2. Preparatory theorems}

\ssbegin{2.1}{Theorem} Let $\dim V=(n, m)$, $\dim U=(k, l)$, $k\ge n$,
$l\ge m$.  Then the algebra $S^{\bcdot}(U\otimes V)$ considered as a
$\fgl(U)\oplus \fgl (V)$-module can be represented in the form
$$
S^{\bcdot}(U\otimes V)=\mathop{\bigoplus} \limits_{\lambda }U^{\lambda
}\otimes V^{\lambda },
$$
where $\lambda $ runs over the set of Young tableaux such that
$\lambda _{n+1}\le m$ and $U^{\lambda }$ and $V^{\lambda }$ are
irreducible $\fgl (U)$- and $\fgl (V)$-modules corresponding to the
tableaux $\lambda $.
\end{Theorem}

\begin{proof} By \cite{S2} we have the following decompositions:
$$
V^{*\otimes k}=\mathop{\oplus}\limits_{\lambda }( V^{*\lambda }\otimes
S^{\lambda }), \quad U^{\otimes k}=\mathop{\oplus}\limits_{\mu }
(U^{\mu }\otimes S^{\mu }).
$$
Here $S^{\lambda }$ and $S^{\mu }$ are irreducible $\fS_{k}$-modules
corresponding to the tableaux $\lambda$ and $\mu$, respectively.

The following isomorphisms take place:
\begin{align*}
S^{k}(U\otimes V)&=S^{k}(U\otimes (V^{*})^{*})=S^{k} (\Hom (V^{*},
U))= \cr
&\Hom _{\fS_{k}}(V^{*\otimes k}, U^{\otimes k})= \Hom
_{\fS_{k}}(\mathop{\oplus}\limits_{\lambda }(V^{*\lambda }\otimes
S^{\lambda }), \mathop{\oplus}\limits_{\mu } (U^{\mu }\otimes
S^{\mu}))=\cr
&\mathop{\oplus}\limits_{\lambda, \;\mu }\Hom (V^{*\lambda
}, U^{\mu })\otimes \Hom _{\fS_{k}}(S^{\lambda }, S^{\mu })=
\mathop{\oplus}\limits_{\lambda }U^{\lambda }\otimes (V^{*\lambda
})^{*}= \cr
&\mathop{\oplus}\limits_{\lambda }(U^{\lambda }\otimes
V^{\lambda })
\end{align*}
All these isomorphisms are $\fgl (U)$ and $\fgl (V)$-isomorphisms.
\end{proof}

\ssec{2.2} The following theorem is similar to Theorem
II.5.A from \cite{Wy}.

\begin{Theorem} Let $\fg $ be a Lie subsuperalgebra in $\fgl (V)$.  If
$\fM$ is a basic system of $\fg$-invariants in $\fA^{n, m}_{n, m}$, then
$\fM$ is also a basic system of invariants in the algebra $\fA^{p,
q}_{k, l}$ for any $k, p\ge n$ and $ q, l\ge m$.
\end{Theorem}

\begin{proof} Let $U_{1}\subset U$, $W_{1}\subset W$ and $\dim U_{1}=
\dim W_{1}=(n, m)$.  By Theorem 2.1 we have
$$
S^{\bcdot}(U_{1}\otimes V\oplus V^{*}\otimes
W_{1})=S^{\bcdot}(U_{1}\otimes V) \otimes S^{\bcdot}(V^{*}\otimes
W_{1})=\mathop{\oplus}\limits_{\lambda, \; \mu } (U^{\lambda
}_{1}\otimes V^{\lambda }\otimes V^{*\mu }\otimes W^{\mu }_{1}),
$$
where $\lambda , \mu $ are Young tableaux
such that $\lambda _{n+1}\le m$, $\mu _{n+1}\le m$.

Similarly,
$$
S^{\bcdot}(U\otimes V\oplus V^{*}\otimes W)=\mathop{\oplus}\limits
_{\lambda, \; \mu } (U^{\lambda }\otimes V^{\lambda }\otimes V^{*\mu }
\otimes W^{\mu }),
$$
where $\lambda$ and $\mu $ are the same here as in the above
expansion.

The embeddings $U_{1}\hookrightarrow U$ and $W_{1}\hookrightarrow W$
induce an embedding $\varphi : \fA^{n, m}_{n, m}\hookrightarrow
\fA^{p, q}_{k, l}$ which is a $\fgl (U_{1})\oplus \fgl
(W_{1})$-homomorphism.  Set
$$
(\fA^{p, q}_{k, l})_{\lambda , \mu }=U^{\lambda }\otimes
V^{\lambda }\otimes V^{*\mu }\otimes W^{\mu }.
$$
So $\varphi ((\fA^{n, m}_{n, m})_{\lambda , \mu }) \subset (\fA^{p,
q}_{k, l})_{\lambda , \mu }$ and since $\varphi $ is also a $\fgl
(V)$-homomorphism, then $\varphi ((\fA^{n, m}_{n, m})^\fg )\subset
(\fA^{p, q}_{k, l})^\fg $.  The elements of the space $(\fA^{p, q}_{k,
l})^\fg _{\lambda , \mu }$ will be called $\fg$-{\it invariants of
type} $(\lambda , \mu )$.

Let $f$ be a $\fg$-invariant highest vector of type $(\lambda , \mu )$
with respect to a Borel subalgebra in $\fgl (U_{1})\oplus \fgl
(W_{1})$.  If we select a Borel subalgebra in $\fgl (U)\oplus \fgl
(W)$ that preserves $U_{1}$ and $W_{1}$, then $f$ is still a highest
vector for such a subalgebra.  This proves that as a $\fgl (U)\oplus
\fgl (W)$-module $(\fA^{k, l}_{p, q})^\fg _{\lambda , \mu }$ is
generated by the subspace $\varphi ((\fA^{n, m}_{n, m})^\fg _{\lambda
, \mu })$ and the theorem is proved.
\end{proof}

\begin{Remark} One can similarly show that

a) if $\fM$ is a basic system of invariants for $\fA^{n, m}$,
then it is a basic system of invariants for $\fA^{p, q}$, where
$p\ge n$, $ q\ge m$;

b) if $\fM$ is a basic system of invariants for $\fA^{n}_{n}$, then
it is also a basic system for any algebra $\fA^{p}_{k}$, where
$p, k\ge n$.
\end{Remark}

Let $A$ be a supercommutative superalgebra over $\Cee$, $L$ a
$\fg$-module and $L_{A}=(L\otimes A)_\ev$, $\fg _{A}= (\fg \otimes
A)_\ev$.  Then the elements of $S^{\bcdot}(L^{*})$ can be considered as
functions on $L_{A}$ with values in $A$.  Let $l\in L_{A}=(L\otimes
A)_\ev=(\Hom (L^{*}, A))_\ev$.  Hence, $l$ determines a homomorphism
$\varphi _{l} : S^{\bcdot}(L^{*})\tto A$.  For $f\in S^{\bcdot}(L^{*})$ set
$f(l)=\varphi _{l}(f)$.  Notice that $\fg _{A}$ naturally acts on
$L_{A}$ and on the algebra of functions on $L_{A}$.

\ssec{2.3.  How to describe $\fg$-invariants in terms of the point
functor} The following result from \cite{S1} essentially means that if $A$
is a Grassmann superalgebra with sufficiently large number of
generators, then for the description of invariants of $\fg $ or the
corresponding Lie supergroup $G$ it suffices to confine ourselves to
$A$-points. We recall the language of points in Appendix.

\begin{Theorem} Let $A$ be a Grassmann superalgebra with the number of
generators greater than $\dim L_{\od}$.  An element of
$S^{\bcdot}(L^{*})$ is a $\fg$-invariant if and only if, considered as
a function on $L_{A}$, it is invariant with respect to $\fg _{A}$.
\end{Theorem}

\begin{Corollary} Let $G_{A}$ be the connected Lie group corresponding
to the Lie algebra $\fg _{A}$.  Then an element of $S^{\bcdot}(L^{*})$ is a
$\fg$-invariant if and only if as a function on $L_{A}$ this element
is $G_{A}$-invariant.
\end{Corollary}

Let
$$
L=V^{p}\oplus \Pi (V)^{q}\oplus V^{*k}\oplus \Pi (V)^{*l}.
$$
Then $S^{\bcdot}(L^{*})=\fA^{p, q}_{k, l}$ and the space $L_{A}$
can be considered as the set of collections
$$
\fL=(v_{1}, \dots , v_{p}, v_{\od}, \dots , v_{\bar q},
v^{*}_{1}, \dots , v^{*}_{k}, v^{*}_{\od}, \dots ,
v^{*}_{\bar l}),
$$
where $v_{s}\in V\otimes A, \; v^{*}_{t}\in \Hom _{A} (V\otimes A, A)$
and their parities coincide with parities of the corresponding
indices.  The vectors will be expressed by means of right coordinates
and covectors by means of left coordinates:
$$
v_{s}=\sum e_{i}a^{*}_{is}, \quad v^{*}_{t}=\sum a_{ti}e^{*}_{i}.
$$
If we consider the elements of the algebra $\fA^{p, q}_{k, l}$
as functions depending on $\fL$, then
$$
x^{*}_{is}(\fL)=a^{*}_{is}, \quad
x_{ti}(\fL)=a_{ti}.
$$

In $\fgl(V)$, introduce a $\Zee$-grading by setting
\begin{align*}
\fgl (V)_{+}&=\{\cA \in \fgl (V):\;\cA V_\ev=0,
\cA V_{\od}\subset V_\ev\}, \cr
\fgl (V)_{0}&=\fgl (V)_\ev, \cr
\fgl (V)_{-}&=\{\cA \in \fgl (V):\;\cA V_{\od}=0,
\cA V_\ev\subset V_{\od}\}.
\end{align*}
Denote by $\fb _{+}(V)$ the Borel subalgebra which consists
of even upper triangular matrices in the basis $\{e_{i}\}_{i\in I}$
and let $\fb _{-}(V)$ be the set of even lower triangular matrices.

We will apply similar notations to $\fgl (U)$ and $\fgl (W)$.

\section*{\S 3. Invariants for the Lie superalgebra $\fgl (V)$}

{\bf Prooof of  Theorem 1.1}.  By Theorem 2.2 it suffices to consider the
case of the algebra $\fA^{n, m}_{n, m}$.  By Corollary to Theorem 2.3
we have to consider functions on collections
$$
\fL=(v_{1}, \dots , v_{n}, v_{\od}, \dots , v_{\overline m},
v^{*}_{1}, \dots , v^{*}_{n}, v^{*}_{\od}, \dots ,
v^{*}_{\overline m})
$$
contained in the algebra generated by coordinate
functions and invariant with respect to the Lie group $GL
(V\otimes A)$. Denote by $M$ the set of collections such that the
vectors
$$
(v_{1}, \dots , v_{n}, v_{\od}, \dots , v_{\overline m})
$$
form a basis in $V\otimes A$. If we consider $M$ as an algebraic
variety, then it is dense in the space of all collections in Zariski
topology. If $f$ is a $GL(V\otimes A)$-invariant and $\fL\in M$, then
there exists $g\in GL (V\otimes A)$ such that $gv_{i}=e_{i}$ for $ i\in
I$; therefore,
$$
f(\fL)=f(g\fL)=f(e_{1}, \dots , e_{\overline m}, gv^{*}_{1}, \dots ,
gv^{*}_{\overline m})
$$
and $f(\fL)$ is a polynomial in coordinates of the vectors $gv^{*}_{t}$.
But
$$
(gv^{*}_{t}, e_{i})=(v^{*}_{t}, g^{-1}e_{i})=(v^{*}_{t}, v_{i}).
$$
Theorem is proved. Its corollary is true since the
polarization operators turn innner products into inner
products. \qed

\section*{\S 4. Invariants of the Lie superalgebra $\fsl (V)$}

\ssec{4.1} By the same reasons as for $\fgl (V)$ it suffices to
confine ourselves to the case of the algebra $\fA^{n, m}_{n, m}$.
First, let us find out for which $\lambda , \mu $ there exist
invariants in $(\fA^{n, m}_{n, m})_{\lambda , \mu }$ and then let us
construct an invariant of type $(\lambda , \mu )$. The tableaux $\lambda
$ and $\mu $ are called {\it equivalent} if the modules $V^{\lambda }$
and $V^{\mu }$ have the same highest weight as $\fsl (V)$-modules.

\begin{Lemma} The tableaux $\lambda $ and $\mu$ are equivalent if and
only if one of the following two case takes place:

{\em a)} $\lambda =\mu$;

{\em b)} $\lambda \neq \mu $ and both $\lambda $ and $\mu $ contain a
rectangle of size $n\times m$ such that there exists $k\in \Zee_+$
such that if we delete $k$ cells from the first $m$ columns of
$\lambda$ and add these $k$ cells to each of the first $n$ rows of
$\lambda$ we get $\mu$.  If $k<0$, then we delete the cells from the
rows and add them to the columns.
\end{Lemma}

\begin{proof} The case a) is obvious.

b) Let $\lambda \neq \mu $ and $\chi _{\lambda }, \chi _{\mu }$ be
highest weights of modules $V^{\lambda }, V^{\mu }$ with respect to
$\fb _{+}(V)$ and let us take the coordinates of the highest weight
with respect to the Cartan subalgebra consisting of diagonal matrix
units.

If $\gamma =(1, \dots , 1, -1, \dots , -1)$, then $\chi _{\lambda
}-\chi _{\mu }=k\gamma $ where $k\in \Zee$.  Let $k>0$; then $(\chi
_{\mu })_{\overline m}>0$ and $\mu _{n}\ge m$.  Then $\lambda _{n}=k+\mu
_{n}>m$, i.e., both tableaux contain an $n\times m$ rectangle.  The
case $k<0$ is treated similarly.  Now, the statement of Lemma is
completely proved.
\end{proof}

\ssbegin{4.2}{Lemma} $\dim (V^{\lambda }\otimes
V^{*\mu })^{\fsl (V)}=\left\{\begin{matrix}
1&\text{ if $\lambda$ and $\mu$ are
equivalent}\cr
0 &\text{otherwise.}\cr\end{matrix}\right .$
\end{Lemma}

\begin{proof}
$$
(V^{\lambda }\otimes V^{*\mu })^{\fsl (V)}=
\Hom _{\fsl (V)}(V^{\mu }, V^{\lambda }).
$$
Thus, $\fsl (V)$-invariants of type $(\lambda , \mu )$ distinct from
$\fgl (V)$-invariants only exist if $\lambda $ and $\mu $ correspond to
typical modules and are equivalent.
\end{proof}

\ssbegin{4.3}{Lemma} Let $M$ and $N$ be finite dimensional
$\fgl (V)_{0}$-modules. Set $\fgl (V)_{+}M=0$, $\fgl (V)_{-}N=0$. Then
$$
\ind^{\fgl (V)}_{\fgl (V)_{0}\oplus
\fgl (V)_{+}}(M)\otimes \ind^{\fgl (V)}_{\fgl (V)_{0}
\oplus \fgl (V)_{-}}(N)=
\ind^{\fgl (V)}_{\fgl (V)_{0}}(M\otimes N). \eqno{(*)}
$$
\end{Lemma}

\begin{proof} Since the dimensions of both modules in $(*)$ are
the same, it suffices to show that the submodule in the
left-hand side of the above equality generated by $M\otimes N$
coincides with the whole module.

Select bases $\{X_{\alpha }\}_{\alpha >0}$ in $\fgl (V)_{+}$ and
$\{Y_{\beta }\}_{\beta <0}$ in $\fgl (V)_{-}$.  Let $L$ be the $\fgl
(V)$-submodule generated by $M\otimes N$.  Consider an element
$$
u=Y_{\beta _{1}}\dots
Y_{\beta _{l}}m\otimes X_{\alpha _{1}}\dots
X_{\alpha _{k}}n, \; \text{ where }\; m\in M, \; n\in N.
$$
Let us prove by induction on $k+l$ that $u\in L$. For $k+l=0$ the
statement is obvious. Let $k+l>0$ and
$$
\tilde u=Y_{\beta _{2}}\dots
Y_{\beta _{l}}m\otimes X_{\alpha _{1}}\dots
X_{\alpha _{k}}n.
$$
By inductive hypothesis $\tilde u\in L$; hence,
$Y_{\beta _{1}}\tilde u\in L$. Furthermore,
$$
u=Y_{\beta _{1}}\tilde u\pm Y_{\beta _{2}}\dots
Y_{\beta _{l}}m\otimes Y_{\beta _{1}}X_{\alpha _{1}}\dots
X_{\alpha _{k}}n
$$
and
$$
\begin{gathered}
Y_{\beta _{1}}X_{\alpha _{1}}\dots
X_{\alpha _{k}}n= \cr
[Y_{\beta _{1}}, X_{\alpha _{1}}]X_{\alpha _{2}}\dots
X_{\alpha _{k}}n-X_{\alpha _{1}}[Y_{\beta _{1}},
X_{\alpha _{2}}]X_{\alpha _{3}}\dots
X_{\alpha _{k}}n+\dots \cr
\pm X_{\alpha _{1}}\dots
X_{\alpha _{k-1}}[Y_{\beta _{1}}, X_{\alpha _{k}}]n.
\end{gathered}
$$
By induction we have $Y_{\beta _{2}}\dots
Y_{\beta _{l}}m\otimes Y_{\beta _{1}}X_{\alpha _{1}}\dots
X_{\alpha _{k}}n\in L$.
\end{proof}

\ssbegin{4.4}{Lemma} Let $\fg =\fgl(V)$ or $\fsl (V)$ and $L$ be a
$\fg$-module.  If $u\in L$ is a $\fg _{0}$-invariant, then
$$
\prod _{\alpha }X_{\alpha }\prod _{\beta }X_{\beta }u(
\prod _{\beta }Y_{\beta }\prod _{\alpha }X_{\alpha }u)
$$
is a $\fg$-invariant (perhaps equal to zero).
\end{Lemma}

{\bf Proof:} the straightforward verification with the
help of the multiplication table for $\fg$. \qed

Let $\lambda =(m^{n+k}), \; \mu =(m+k)^{n}$; then by Lemma 4.2 in the
algebra $\fA^{n, 0}_{0, m}=S^{\bcdot}(V^{*n}\oplus \Pi (V)^{m})$ there
exists an invariant of type $(\lambda , \mu )$.  It is not difficult
to see that this invariant is unique up to a constant factor.  The
submodule generated by this invariant is isomorphic to $(\Ber~ \Pi
(V))^{\otimes k}$.

Similarly, in the algebra $\fA^{0, m}_{n, 0}=S^{\bcdot}(V^{n}\oplus \Pi
(V)^{*m})$ there exists a unique $\fsl (V)$-invariant of type $(\mu ,
\lambda )$ and the module generated by it is isomorphic to $(\Ber~
V)^{\otimes k}$.

An explicit description of these invariants is given in Theorem 1.2
(polynomials $f_{k}$).

\ssbegin{4.5}{Lemma} For $k\in \Nee$ the polynomials
$$
f_{k}=(\Delta ^{*})^{k}\omega ^{k}
\prod \limits_{t\in I_{\od}, \;
s\in I_\ev}(v^{*}_{t}, v_{s}), \quad \text{and}\quad
f_{-k}=\Delta ^{k}(\omega ^{*})^{k}
\prod \limits_{t\in I_\ev, \;
s\in I_{\od}}(v^{*}_{t}, v_{s})
$$
are $\fsl (V)$-invariant.
\end{Lemma}

\begin{proof} Consider $f_k$ for $k>0$ (the case $k<0$ is similar).
Select bases $\{X_{\alpha }\}_{\alpha >0}$ in $\fgl (V)_{+}$ and
$\{Y_{\beta }\}_{\beta <0}$ in $\fgl (V)_{-}$ and denote
$X=\mathop{\prod}\limits_{\alpha }X_{\alpha }$,
$Y=\mathop{\prod}\limits_{\beta }Y_{\beta }$.  Further on, introduce
polynomials
$$
\Pi ^{*}_{10}=\mathop{\prod}\limits _{i\in I_{\od}, \;
s\in I_\ev} x^{*}_{is}, \quad \Pi _{10}=
\mathop{\prod}\limits _{t\in I_{\od}, \;
i\in I_\ev} x_{ti}.
$$
We have
$$
X(\Delta ^{*m+k})=c_{1}\Delta ^{*k}\Pi ^{*}_{10};
\eqno{(4.5.1)}
$$
$$
YX(\Delta ^{*m}\Pi _{10})=c_{2}
\mathop{\prod}\limits_{t\in I_{\od}, \;
s\in I_\ev} (v^{*}_{t}, v_{s});
\eqno{(4.5.2)}
$$
$$
YX(\Delta ^{*m+k}\omega ^{k}\Pi _{10})=c_{3}\Delta ^{*k}\omega ^{k}
\prod\limits_{t\in I_{\od}, \; s\in I_\ev}(v^{*}_{t}, v_{s}),
\eqno{(4.5.3)}
$$
where $c_{1}, c_{2}, c_{3}$ are nonzero constants.

Indeed, consider $S^{\bcdot}(V^{*n})_{\mu }$, where $\mu =(m+k)^{n}$.
The elements $\Delta ^{*m+k}$ and $\Delta ^{*k}\Pi ^{*}_{10}$ belong
to $S^{\bcdot}(V^{*n})_{\mu }$ and
$$
\fgl (V)_-(\Delta ^{*m+k})=\fgl (V)_+(\Delta ^{*k}\Pi ^{*}_{10})=0.
$$
Since the module corresponding to $\mu $ is a typical one,
the equality (4.5.1) holds.

Now, consider $S^{\bcdot}(V^{*n}\oplus \Pi (V)^{m})_{\lambda , \lambda }$,
where $\lambda =(m^{n})$; we have $\dim U^{\lambda }=\dim
W^{\lambda}=1$ for this $\lambda $ and, therefore, there exists only
one invariant of type $(\lambda , \lambda )$, namely,
$$
\prod\limits_{t\in I_{\od}, \; s\in I_\ev} (v^{*}_{t}, v_{s}).
$$
On the other hand, $\Delta ^{*m}\Pi _{10}$ belongs to
$V^{\lambda }\otimes V^{*\lambda }$ and is a
$\fgl (V)_{0}$-invariant. Lemmas 4.3 and 4.4 imply (4.5.2).

Finally:
\begin{align*}
&YX(\Delta ^{*m+k}\omega _{k}\Pi _{10})=Y(c_{1}\Delta ^{*k}
\Pi ^{*}_{10}\omega ^{k}\Pi _{10})= \cr
&c_{1}\Delta ^{*k}\omega ^{k}Y(\Pi ^{*}_{10}\Pi _{10})=
c_{1}c_{2}\Delta ^{*k}\omega ^{k}\prod
\limits_{t\in I_{\od}, \; s\in I_\ev} (v^{*}_{t}, v_{s})
\end{align*}
and by Lemma 4.4 this expression is a $\fsl(V)$-invariant.
Proof is completed.
\end{proof}

\ssec{4.6.  Proof of Theorem 1.2} It suffices to construct an
invariant polynomial in $(\fA^{n, m}_{n, m})_{\lambda , \mu }$ for
$\lambda $ and $\mu $ indicated in Lemma 4.1 that only depends on the
polynomials $f_{k}$ and the inner products.  Let $k>0$ and $\lambda $
and $\mu $ be chosen as on Fig. 4.6.

\centerline{\epsfbox{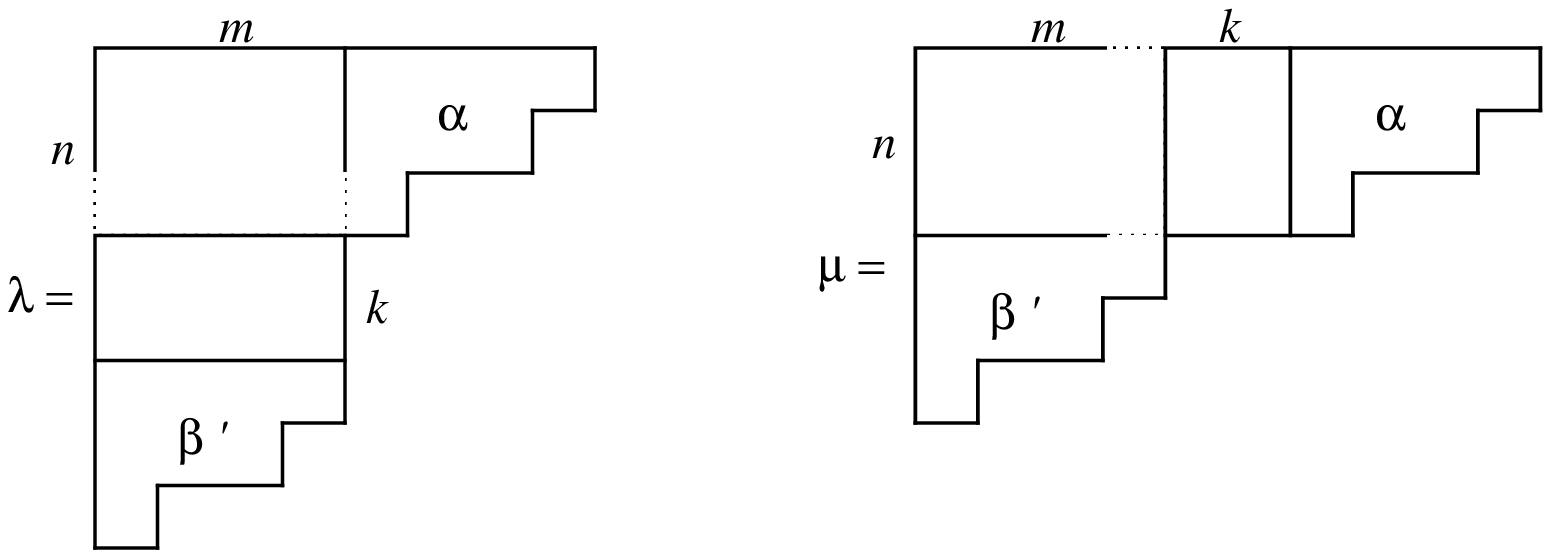}}
\centerline{Figure 4.6}

Let $\varphi _{\alpha }$ be a $\fb _{+}(U)\oplus \fb _{+}(W)$-highest
invariant of type $(\alpha , \alpha )$ in the algebra
$S^{\bcdot}(V^{n}\oplus
V^{*n})$ and $\psi _{\beta }$ be a similar invariant of type $(\beta ,
\beta )$ in the algebra $S^{\bcdot}(\Pi (V)^{m} \oplus \Pi (V)^{*m})$.
Then we
can verify that $f_{k}\varphi _{\alpha }\psi _{\beta }$ is the $(\fb
_{+}(U)\oplus \fgl (U)_{+})\oplus (\fb _{+}(W) \oplus \fgl
(W)_{-})$-highest vector.  Clearly, it is an $\fsl (V)$-invariant and
its weight corresponds to the tableaux $\lambda $ and $\mu $.

The case $k<0$ is treated similarly.

\section*{\S 5. Invariants of the Lie superalgebra $\fosp (V)$}

\ssec{5.0} Let $\dim V=(n, 2r)$ and $\fosp (V)$ be the Lie
superalgebra described in \S 1.  By Remark 2.2 it suffices to confine
ourselves to the algebra $\fA^{n, 2r}$.  From the point of view of
Theorem 2.3 we have to describe the polynomials that depend on the set
$$
(v_{1}, \dots , v_{n}, v_{\od}, \dots , v_{\overline{2r}})
$$
and invariant with respect to the simply connected Lie group $G_{A}$
whose Lie algebra is $(\fosp (V)\otimes A)_\ev$, where $A$ is the
Grassmann superalgebra with a sufficiently large number of generators.
Denote by $\OSp(V\otimes A)$ the subgroup of $\GL(V\otimes A)$ whose
elements preserve the inner product
$$
(v_{s}, v_{t})=\sum^{n}_{i=1}x^{*}_{is}x^{*}_{n-i+1, t}+(-1)^{p(s)}
\sum^{2r}_{j=1}(x^{*}_{\overline{m-j+1}, s}x^{*}_{\bar j, t}
-x^{*}_{\bar j, s}x^{*}_{\overline{m-j+1}, t})
$$
and by $S\OSp(V\otimes A)$ the subgroup of $\OSp(V\otimes A)$
consisting of transformations with Berezinian 1.

It is not difficult to verify that $S\OSp(V\otimes A)$ is
precisely the group $G_{A}$ spoken about above.

Denote by $O(V_\ev)$ the orthogonal group that preserves
the form
$$
\sum^{n}_{i=1}x^{*}_{i}x^{*}_{n-i+1}.
$$
It is not difficult to verify that the invariance of an element of
$\fA^{n, m}$ with respect to $\OSp(V\otimes A)$ is equivalent to the
simultaneous invariance with respect to $\OSp(V)$ and $O(V_\ev)$.

First, let us prove several lemmas.

\ssbegin{5.1}{Lemma} Let $M$ be a $\fg (V)_\ev$-module.  Set
$\fgl(V)_{+}(M)=0$.  Then we have an isomorphism of $\fosp (V)$ and
$O(V_\ev)$-bimodules:
$$
\ind^{\fgl (V)}_{\fgl (V)_{0}\oplus \fgl (V)_{+}}(M)\simeq
\ind^{\fosp (V)}_{\fosp (V)_\ev}(M).
$$
\end{Lemma}

\begin{proof} Let us describe a basis in $\fosp (V)_{1}$. A
nondegenerate form determines a map in $\fgl (V)$:
$$
\cA \mapsto \overline{\cA } \; \; \text {and ~
$\overline{\overline{\cA }}=(-1)^{p(\cA )}\cA $.}
$$
If $S$ is a matrix of the form and $P$ is the matrix of $\cA $ in the
same basis then
$$
\overline P=S^{-1}P^{st}S, \; \text{where $P^{st}$ is the
supertransposed matrix.}
$$
If $\{X_{\alpha }\}$ is a basis in $\fgl (V)_{-}$,
then $\{X_{\alpha }-\overline X_{\alpha }\}$ is a basis in
$\fosp (V)_{\od}$ such that $\overline X_{\alpha }\in
\fgl (V)_{+}$. Let
$$
\varphi : \ind^{\fosp (V)}_{\fosp (V)_{0}}(M)
\tto\ind^{\fgl (V)}_{\fgl (V)_{0}\oplus
\fgl (V)_{+}}=L
$$
be a homomorphism induced by the natural embedding $M\hookrightarrow
L$.  On $L$, there exists a filtration $L_{0}\subset L_{1}\subset
\dots \subset L_{N}$, where $L_{k}$ is the linear hull of $f(X_{\alpha
})m$ and where $m\in M$, $\deg f\le k$.

Let us prove by induction that $L_{k}\subset {\rm Im} \varphi $.

$k=0$ is an obvious case.  Let $L_{k}\subset {\rm Im} \varphi $, then
$(X_{\alpha }-\overline X_{\alpha })L_{k}\subset {\rm Im} \varphi $ but
$\overline X_{\alpha }L_{k}\subset L_{k-1}$.  Therefore, $X_{\alpha
}L_{k}\subset {\rm Im} \varphi $ or $L_{k+1}\subset {\rm Im} \varphi
$.

The statement on $O(V_\ev)$-modules is obvious.  Lemma is proved.
\end{proof}

\ssbegin{5.2}{Lemma} Let $\fg $ be a finite dimensional Lie
superalgebra and the representation of $\fg _{\ev}$ in $\Lambda ^{\dim
\fg_{\od}}(\fg _{\od})$ be trivial.  Then for a finite-dimensional
$\fg_{\ev}$-module $M$ there exists an isomorphism of vector spaces:
$$
\left(\ind^{\fg }_{\fg _\ev}
(M)\right)^{\fg }\simeq M^{\fg _\ev}.
$$
\end{Lemma}

\begin{proof} Let us prove that $(\ind^{\fg }_{\fg _\ev}
(M))^{*}_{\fg }\simeq\ind^{\fg }_{\fg _\ev} (M^{*})$ as $\fg$-modules.
Let $\fg _{\od}=\Span(\xi _{1}, \dots , \xi _{p})$.  Then $L=\ind^{\fg
}_{\fg _\ev}(M)$ has a natural filtration with $\fg _\ev$-modules and
as in the above lemma $L_{0}\subset L_{1}\subset \dots \subset
L_{p}=L$.

The map
$$
M\tto L_{p}/L_{p-1},\;\; m\mapsto \xi _{1}\dots
\xi _{p}m\; ({\rm mod}\; L_{p-1})
$$
induces an isomorphism of $\fg _{\ev}$-modules: $M^{*}
(L_{p}/L_{p-1})^{*}$. Therefore, we have an embedding of $\fg
_\ev$-modules
$$
M^{*}\tto (L_{p}/L_{p-1})^{*}\tto L^{*}.
$$
This map induces a homomorphism $\ind^{\fg }_{ \fg
_\ev}(M^{*})\tto L^{*}$.  Consider the $\fg$-invariant
bilinear form corresponding to this homomorphism:
\begin{align*}
&\ind^{\fg }_{\fg _\ev}(M^{*})\times \ind^{
\fg }_{\fg _\ev}(M)\tto\Cee \cr
&(m^{*}, \xi _{1}\dots
\xi _{p}m)=m^{*}(m) \; \text{ for }\; m^{*}\in M^{*}, \; m\in M.
\end{align*}
Let $u$ be a nonzero element from the left kernel of the form.  There
exists a filtration on the module $T=\ind^{\fg }_{\fg _{\ev}}(M^{*})$;
same as on $L$.  Let $u\in T_{k}$ but $u\not\in T_{k-1}$.  Then
$$
u=\sum \xi _{i_{1}}\dots \xi _{i_{k}}m^{*}_{i_{1}\dots i_{k}}+u_{k-1}.
$$
Set
$$
v=\xi _{j_{1}}, \dots , \xi _{j_{l}}m, \quad \text{where}\quad
\{j_{1}\dots j_{l}\}=[1, \dots , p]\setminus \{i_{1}, \dots , i_{k}\}.
$$
Then
$$
(u, v)=(\xi _{i_{1}}\dots \xi _{i_{k}}m^{*}_{i_{1}\dots i_{k}}, v)=\pm
m^{*}_{i_{1}\dots i_{k}}(m)=0.
$$
Since $m$ is arbitrary, $m^{*}_{i_{1}\dots i_{k}}=0$ and, therefore,
$u=u_{k-1}\in T_{k-1}$; contradiction. Hence, $u=0$. Further,
$$
(\ind^{\fg }_{\fg _{\ev}}(M))^{\fg }=(\ind^{
\fg }_{\fg _\ev}(M^{*}))^{*\fg }=(M^{*})^{*\fg _\ev}=M^{\fg_\ev}.
$$
\end{proof}

\begin{Remark} 1) If $\fg _{\ev}\supset \fo(n)$ and $M$ such that
$\fg_{\od}$ are $O(n)$-modules and the $O(n)$-action in $\Lambda
^{p}(\fg_{\od})$ is trivial, then the statement of the lemma remains
valid for the mutual $\fg$- and $O(n)$-invariants and mutual
$\fg_\ev$- and $O(n)$-invariants.

2) The following refinement of the lemma can be obtained: {\sl if
$m\in M$ is a $\fg_{\ev}$-invariant, then the corresponding
$\fg$-invariant vector $u$ is of the form}
$$
u=\xi _{1}\dots \xi _{p}m+u_{p-1},\; \text{ where }\; u_{p-1}\in L_{p-1}.
$$
\end{Remark}

\ssec{5.3. Description of invariants} First, consider
$\OSp(V\otimes A)$-invariants.

\begin{Theorem} Any $\OSp(V\otimes A)$-invariant element from $\fA^{p,
q}$ is a polynomial in inner products $(v_{s}, v_{t})$, where $s, t\in S$.
\end{Theorem}

\begin{proof} Induction on $\dim V_\ev$.

If $\dim V_\ev=0$, then Theorem is proved in \cite{Wy}.  Let $\dim
V_\ev=n>0$.  It suffices to show that any invariant of type $\lambda $
with $\lambda _{n+1}\le 2r$ can be expressed in terms of inner
products.  First, let $\lambda$ satisfy the condition $\lambda_{n}\le 2r$.

Consider the algebra
$$
\fA^{n-1, 2r}=\mathop{\oplus}\limits_{\lambda _{n}\le 2r}
V^{*\lambda }\otimes W^{\lambda }.
$$
If $\{v_{1}, \dots , v_{n-1}, v_{\od}, \dots , v_{\overline{2r}}\}$ is a
collection of vectors in general position from $V\otimes A$, then
after an orthogonalization we may assume that there exists $g\in
OSp(V\otimes A)$ such that
$$
g\Span(v_{1}, \dots , v_{n-1}, v_{\od}, \dots ,
v_{\overline{2r}})=H=\Span(e_{1}, \dots , e_{n-1}, e_{\od}, \dots ,
e_{\overline{2r}}).
$$
Let $f\in \fA^{n-1, 2r}$ be an invariant with respect to $OSp(V\otimes
A)$ and $\overline f$ the restriction of $f$ to $H$.  By the inductive
hypothesis $\overline f$ is a function in inner products $(\overline v_{s},
\overline
v_{t}), \; \overline v_{s}, \overline v_{t}\in H$.  Hence,
\begin{align*}
&f(v_{1}, \dots , v_{n-1}, v_{\od}, \dots , v_{\overline{2r}})=\overline
f(gv_{1}, \dots , gv_{\overline{2r}})=\cr
&F((gv_{s}, gv_{t}))=F( (v_{s},
v_{t})) \; \text{ for }\; s, t\in \setminus \{n\}.
\end{align*}
Now, let $\lambda _{n}>2r$ but $\lambda _{n+1}\le 2r $.  Then the
$\fgl (V)$-module $V^{*\lambda }$ is a typical one and therefore
$$
V^{*\lambda }=\ind^{\fgl (V)}_{\fgl (V)_{0}\oplus
\fgl (V)_{+}}(M),
$$
where $M$ is an irreducible $\fgl (V)_{0}$-module.  If $\lambda $ is
of the form shown on Fig.  5.3

\centerline{\epsfbox{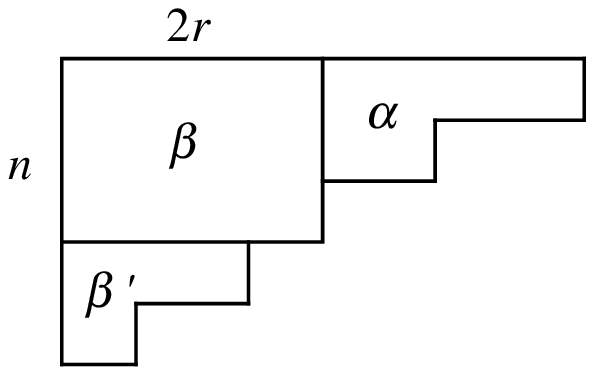}}
\centerline{Figure 5.3}

\noindent then $M=V^{*\alpha }_\ev\otimes
V^{*\beta +\delta '}_{\od}$.

As is not difficult to verify for the orthogonal case and, similarly,
for the symplectic case,
$$
\begin{gathered}
\dim (V^{*\alpha }_{
\bar 0})^{\fo(V_\ev)}=\left\{\begin{matrix}1&\text{ if $\alpha $
is even}\cr
0 &\text{otherwise},\end{matrix}\right .\cr
 \dim (V^{*\beta +\delta '}_{\od})^{\fsp}(V_{\od})=\left\{\begin{matrix}
1&\text{if $(\beta +\delta ')'$ is even}\cr
0 &\text{otherwise.}\end{matrix}\right .
\end{gathered}
$$
These conditions are equivalent to the fact that $\lambda $ is even
(all rows are of an even length). Lemmas 5.1 and 5.2 imply that if
$\lambda $ is typical, then
$$
\dim (V^{*\lambda })^{\OSp(V\otimes A)}=\left\{\begin{matrix}
1&\text{if $\lambda $ is even}\cr
0 &\text{otherwise.}\end{matrix}\right .
$$

Further, for $s, t\in I$ the inner products $(v_{s}, v_{t})$ are
algebraically independent. If we consider the algebra
$\Cee [(v_{s}, v_{t})_{s, t\in I}] $ as a $\fgl (W)$-module,
then
$$
\Cee [(v_{s}, v_{t})_{s, t\in I}]=
\mathop{\oplus}\limits_{\lambda _{n+1}\le 2r} W^{\lambda }.
$$

This is a corollary of a general identity for $\lambda $-rings, see
\cite{M}, \S 5.  This shows that if $\lambda $ is typical and even, then
there exists an invariant of type $\lambda $ depending on inner
products.  The induction is completed and Theorem is proved.
\end{proof}

\ssbegin{5.4}{Lemma} There exists an $\fosp (V)$-invariant $\Omega $
such that
$$
\Omega ^{2}=[\det (v_{s}, v_{t})_{s, t\in I_{0}}]^{2r+1}.
$$
\end{Lemma}

\begin{proof} We have
$$
S^{\bcdot}(V^{*n})=\mathop{\oplus}\limits_{\lambda _{n+1}=0}
V^{*\lambda}\otimes W^{\lambda}.
$$
Let $\lambda =((2r+1)^{n})$; then $\dim W^{\lambda }=1$, the module
$V^{*\lambda }$ is typical, and $V^{*\lambda }=\ind^{\fgl (V)}_{\fgl
(V)_{0}\oplus \fgl (V)_{-}}(M)$, where $\dim M=1$ and $M=\Span(\Delta
^{*}\Pi ^{*}_{10})$.  By Lemmas 5.1 and 5.2 there exists an $\fosp
(V)$-invariant,
$$
\Omega =\prod (X_{\alpha }-\overline X_{\alpha })\Delta ^{*}\Pi ^{*}_{10}+
\tilde \Omega =\prod X_{\alpha }\Delta ^{*}\Pi ^{*}_{10}+\tilde \Omega
_{1}= e\cdot \Delta ^{*2r+1}+\tilde \Omega _{1}.
$$
Therefore, $\Omega ^{2}\neq 0$ but $\Omega ^{2}$ is an $OSp(V\otimes
A)$-invariant and its type is equal to $((4r+2)^{n})$.  However, as is
easy to see, the algebra $S^{\bcdot}(V^{*n})$ has only one $OSp(V\otimes
A)$-invariant of such type, namely, $[\det (v_{s}, v_{t})_{s, t\in
I_{0}}]^{2r+1}$.  Lemma is proved.
\end{proof}

\ssec{5.5.  Proof of Theorem 1.3} Let $f$ be an $\fosp (V)$-invariant
but not an $O(V_\ev)$-invariant.  Let $f$ depend on $n-1$ even and
$2r$ odd vectors.  Let these vectors be in general position.  Then as
in Theorem 5.3 there exists $g\in OSp(V\otimes A)$ such that
$$
g\Span(v_{1}, \dots , v_{\overline{2r}})=\Span(e_{1}, \dots ,
e_{\overline{2r}}).
$$
Let $he_{n}=-e_{n}$ and $he_{i}=e_{i}, \; i\neq n$.  Then $\Ber~
(h)=-1$ and $f(hg\fL)=-f(g\fL)$.  On the other hand, $f(hg\fL)=
f(g\fL)$; therefore, $f=0$.  This means that $\fosp (V)$-invariants
distinct from inner products can be of type $\lambda $ which only
corresponds to a typical module.

Therefore, we can apply Lemmas 5.1 and 5.2.  The same arguments as in
Theorem 5.3 yield that $\dim (V^{*\lambda })^{\fosp (V)}=1$ if
$\lambda $ is typical and its first $n$ rows are of odd length whereas
the remaining rows are of even length and $\dim (V^{*\lambda })^{\fosp
(V)}=0$ otherwise.  (We do not take $\OSp(V\otimes A)$-invariants into
account.)  Let $\lambda$ be as on Fig.  5.5.

\centerline{\epsfbox{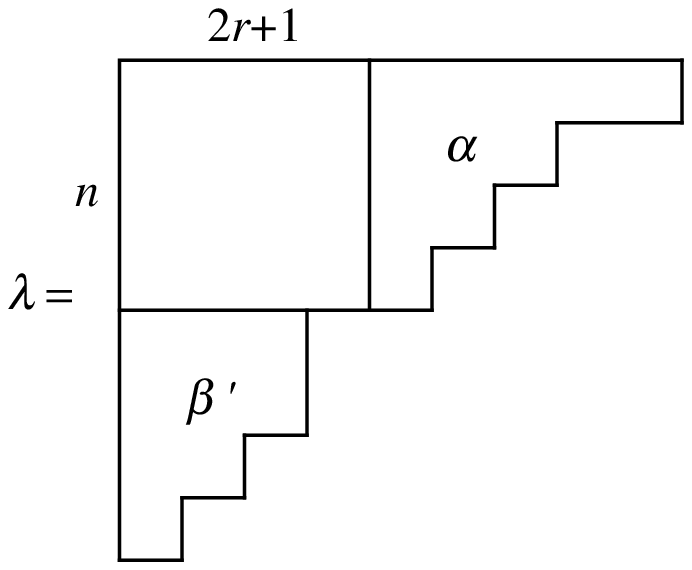}}
\centerline{Figure 5.5}

Let us construct an invariant of type $\lambda $.  Denote by $\varphi
_{\alpha }$ the invariant of type $\alpha $ highest with respect to
$\fb _{+}(W)$ in $S^{\bcdot}(V^{*n})$ and let $\psi _{\beta }$ be the
invariant
of type $\beta $ highest with respect to $\fb _{+}(W)$ in $S^{\bcdot}(\Pi
(V)^{*2r})$.  If $\{D_{j}\}$ is a basis of $\fgl (W)_{-}$, then it is
not difficult to verify that $\psi _{\beta }\prod D_{j}(\Omega \varphi
_{\alpha })$ is $\fb _{+}(W)\oplus \fgl (W)_{-}$-highest and is of
type $\lambda $.  Theorem is proved.
\qed

\section*{\S 6. The invariants of the Lie superalgebra $\fpe(V)$}

\ssec{6.1} Let $\dim V=(n, n)$.  By remark to Theorem 2.2 it suffices
to confine ourselves to the algebra $\fA^{n, n}$.  Denote by
$\Pe(V\otimes A)$ the subgroup in $\GL(V\otimes A)$ whose elements
preserve the inner products
$$
(v_{s}, v_{t})=\sum^{n}_{i=1}((-1)^{p(s)}x^{*}_{is}x^{*}_{ \bar
it}+x^{*}_{\bar is}x^{*}_{it})\; \text{ for }\; s, t\in I_{0}.
$$
This is the group corresponding to the Lie algebra of $A$-points of
the Lie superalgebra $\fpe(V)$.

\begin{Lemma} Let $L$ be an irreducible typical $\fgl (V)$-module,
$$
L=\ind^{\fgl (V)}_{\fgl (V)_{0}\oplus \fgl (V)_{+}}(M)\simeq
\ind^{\fgl (V)}_{\fgl (V)_{0}\oplus \fgl (V)_{-}}(N).
$$
Then there is an isomorphism of vector spaces
$L^{\fsp(V)}=M^{\fsp(V)_\ev}=N^{\fsp(V)_\ev}$.
\end{Lemma}

\begin{proof} The decomposition
$$
\fgl (V)=\fgl (V)_{-}\oplus
\fgl (V)_{0}\oplus \fgl (V)_{+}
$$
induces the decomposition
$$
\fsp(V)=\fsp(V)_{-}\oplus \fsp(V)_{0}
\oplus \fsp(V)_{+}
$$
which is a $\Zee$-grading.  Select a basis $\{Y_{\beta }\}_{1\le \beta
\le n^{2}}$ in $\fgl (V)_{-}$ so that $\{Y_{\beta }\}_{1\le \beta \le
\frac{1}{2}n(n-1)}$ is a basis in $\fsp(V)_{-}$.  Similarly, select a
basis $\{X_{\alpha }\}_{1\le \alpha \le n^{2}}$ in $\fgl (V)_{+}$ so
that $\{X_{\alpha }\}_{1\le \alpha \le \frac{1}{2}n(n+1)}$ is a basis
in $\fsp(V)_{+}$.  Consider two gradings of the module $L$:
$$
\begin{gathered}
L^{+}_{k}=\Span(f(X_{\alpha })n: n\in N, \; \deg f=k\; \text{ for all }\;
\alpha),\cr
L^{-}_{k}=\Span(f(Y_{\beta })m: m\in M, \; \deg f=k\; \text{
for all }\; \beta).
\end{gathered}
$$
It is not difficult to verify that $L^{+}_{k}= L^{-}_{n^{2}-k}$.  Let
$l$ be an $\fsp(V)$-invariant, then $X_{\alpha }l=0$ for $1\le \alpha
\le \frac12n(n+1)$.  Let $X^{+}=\prod X_{\alpha }$ for $\alpha \le
\frac12n(n+1)$; hence, $l=X^{+}f(X_{\alpha })n$ for $n\in N$ and,
therefore, $l=\sum^{}_{r\ge \frac12n(n+1)}l^{+}_{r}$, where
$l^{+}_{r}\in L^{+}_{r}$.

We can similarly verify that $l=\sum\limits^{}_{s\ge
\frac12n(n-1)}l^{-}_{s}$ for $l^{-}_{s}\in L^{-}_{s}$.  This implies
$\sum\limits^{}_{r\ge \frac12n(n+1)}l^{+}_{r}=\sum\limits^{}_{s\ge
\frac12n(n-1)}l^{-}_{s}$.  Since $L^{-}_{s}=L^{+}_{n^{2}-s}$, we get
$l\in L^{+}_{\frac12n(n+1)}=L^{-}_{\frac12n(n-1)}$ and
$l=X^{+}n=X^{-}m$, where $m\in M$ and $X^{-}=\prod Y_{\beta }$ for
$\beta \le \frac12n(n-1)$.

It is clear that $m$ and $n$ are $\fsp(V)_\ev=\fsl
(V_\ev)$-invariants.  Conversely, if $m$ and $n$ are $\fsl
(V_\ev)$-invariants, then the straightforward verification shows that
$X^{+}n$ and $X^{-}m$ are $\fsp(V)$-invariants.  The statement on
bijection is obvious.  The lemma is proved.
\end{proof}

\ssec{6.2. Proof of Theorem 1.4} Induction on $\dim V_\ev$.

The case $n=1$ is straightforward; let theorem hold for $\dim
V_\ev=n-1$.  Let us consider $f\in S^{\bcdot}(V^{*n-1}\oplus \Pi
(V)^{*n-1})^{\Pe(V)}$, where $\dim V=(n, n)$.  For the generic vectors
$$
(v_{1}, \dots , v_{n-1}, v_{\od}, \dots , v_{\overline{n-1}})
$$
there exists $g\in \Pe(V\otimes A)$ such that
$$
gv_{i}\in \Span( e_{1}, \dots , e_{n-1}, e_{\od}, \dots ,
e_{\overline{n-1}})=H.
$$
Let $\overline f$ be the restriction of $f$ to $H$; by the inductive
hypothesis
$f$ is a polynomial in inner products but
$$
f(v_{1}, \dots , v_{\overline{n-1}})=f(gv_{1}, \dots ,
gv_{\overline{n-1}})=F((gv_{s}, gv_{t})_{s, t\in I_{0}})= F((v_{s},
v_{t})_{s, t\in I_{0}}).
$$
Since
$$
S^{\bcdot}(V^{*n-1}\oplus \Pi (V)^{*n-1})=\mathop{\oplus}\limits _{\lambda
_{n}\le n-1} V^{*\lambda }\otimes W^{\lambda }
$$
it remains to demonstrate that the invariants in typical modules can
be expressed in terms of inner products.  Let $\lambda $ be a tableaux
of the form shown on Fig.  6.2 a):

\centerline{\epsfbox{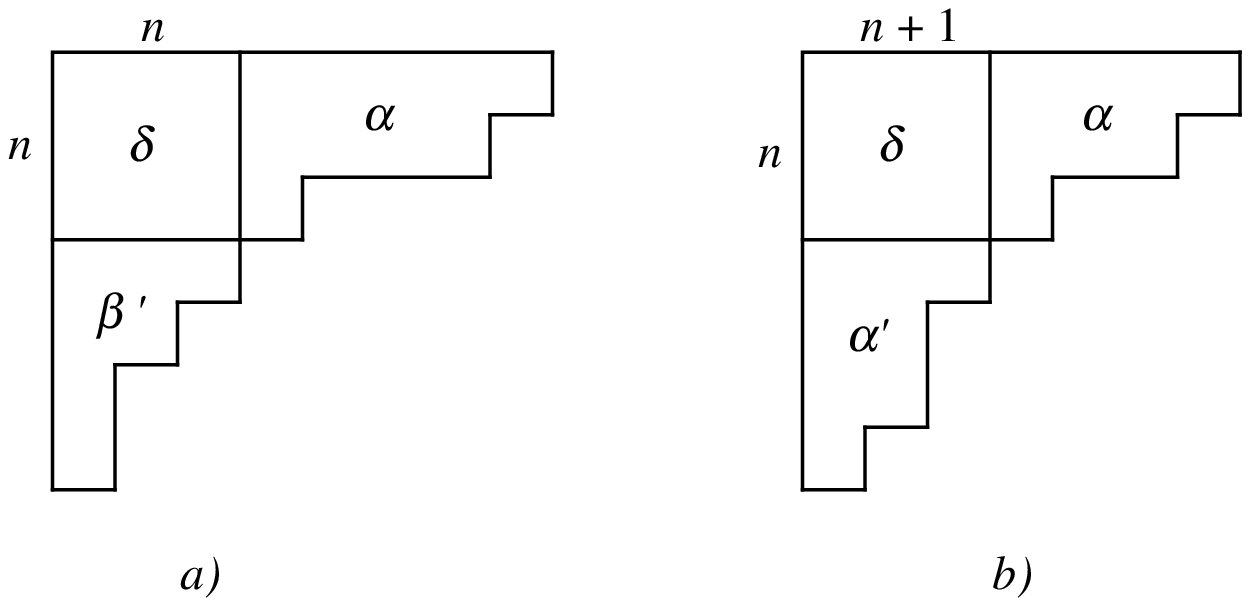}}
\centerline{Figure 6.2 a)\hskip 5 cm Figure  6.2 b)}

\noindent Having applied Lemma 6.1 we see that an invariant of type
$\lambda $ exists if and only if there exists an $\fsl
(V_\ev)$-invariant in $M=V^{*\alpha }_\ev\otimes V^{*\beta +\delta
'}_{\od}$.  We have an isomorphism of $\fsp(V)_\ev$-modules
$$
M\cong V^{*\alpha }_\ev\otimes V^{\beta +\delta '}_{\bar 0}=\Hom
(V^{\alpha }_\ev, V^{\beta + \delta '}_\ev)
$$
and, therefore, $M$ contains an $\fsl (V_{ \bar 0})$-invariant if
$\beta +\delta '-\alpha $ is a multiple of $\gamma =(1, \dots , 1)$.

Once again by Lemma 6.1 the invariant is of the form $X^{-}m$ and,
since we wish it to be $\fpe(V)$-invariant, we need it to be a $\fgl
(V_\ev)$-invariant.  Its weight is equal to
$$
\beta +\delta '-\alpha -(n-1)\gamma =\beta +n\gamma -
\alpha -(n-1)\gamma =\beta -\alpha +\gamma =0,
$$
i.e., $\alpha =\beta +\gamma $ and the tableaux
$\lambda $ should be of the form shown on Fig. 6.2 b).

Let us explicitly indicate an invariant of such
type $\lambda $. The Lie algebra $\fgl
(W_\ev)\oplus \fgl (W_{\od})$ acts on the algebra
$$
\tilde \fA=\Cee [(v_{i}, v_{j})_{i\in I_\ev, \quad
j\in I_{\od}}]
$$
and with respect to this action $\fA= \oplus _{\alpha }W^{\alpha
}_\ev\otimes W^{\alpha }_{\od}$.  Let $\varphi _{\alpha }$ be a vector
from $\tilde \fA$ of type $\alpha $ highest with respect to $\fb
_{+}(W)$.  Then we can verify that $\prod\limits _{1\le i\le j\le
n}(v_{i}, v_{j})\varphi _{\alpha }$ is a $\fb _{+}(W)\oplus \fgl
(W)_{-}$-highest vector of type $\lambda $.  Theorem is proved.  \qed

\section*{\S 7. The invariants of the Lie superalgebra $\fspe(V)$}

\ssec{7.1} First, let us construct certain $\fspe(V)$-invariant
elements in the algebra $\fA^{n, n}$.

\begin{Lemma} The polynomials $\Delta ^{*k}
\prod \limits_{s\le t,\; s, t\in I_\ev}
(v_{s}, v_{t})$ are
$\fspe(V)$-invariant for $k=1$, $2$, $3$, \dots
\end{Lemma}

\begin{proof} Let us consider $\fA=S^{\bcdot}(V^{*n})$ and let $\lambda
=((n+1)^{n})$. Then, in $\fA$, there exists only one invariant of type
$\lambda$:
$$
\Pi ^{+}=\prod \limits_{s\le t,\; s, t\in I_\ev}
(v_{s}, v_{t}).
$$
On the other hand, by Lemma 6.1 we have
$$
X^{-}(\Delta ^{*}\Pi ^{*}_{10})=\Pi ^{+}
$$
(the notation $\Pi ^{*}_{10}$ is defined in Lemma 4.5).

Further, the vector $\Delta ^{*k+1}\Pi ^{*}_{10}$ is a
$\fspe(V)_\ev$-invariant and $\fgl (V)_{+}(\Delta ^{*k+1}\Pi
^{*}_{10})=0$, therefore, by Lemma 6.1 we get an invariant
$X^{-}(\Delta ^{*k+1}\Pi ^{*}_{10})$.  It is not difficult to verify
that $X^{-}(\Delta ^{*})=0$; hence,
$$
X^{-}(\Delta ^{*k+1}\Pi ^{*}_{10})=\Delta ^{*k}X^{-}
(\Delta ^{*}\Pi ^{*}_{10})=\Delta ^{*k}\Pi ^{+}.
$$
Lemma is proved.
\end{proof}

\ssbegin{7.2}{Lemma} The polynomials
$\omega ^{*k}\prod \limits_{s\le t,\; s, t\in I_\ev}(v_{\bar s},
v_{\bar t})$
are $\fspe(V)$-invariant for $k=1, 2, 3, \dots$.
\end{Lemma}

\begin{proof} Let again $\lambda =((n+1)^{n})$.  In $\fA^{n, n}$,
consider the vector
$$
n=\Delta ^{*}\prod\limits_{i\in I_\ev, \;
s\in I_{\od}} X^{*}_{is}=\Delta ^{*}\Pi ^{*}_{01}.
$$
It is not difficult to verify that $\fgl (W)_{-}(\Delta ^{*}\Pi
^{*}_{01})=0$ and $\Delta ^{*}\Pi ^{*}_{01}$ is a
$\fspe(V)_\ev$-invariant highest with respect to $\fb _{+}(W)\oplus
\fgl (W)_{-}$ and of type $\lambda $.  By Lemma 6.1 $X^{+}n$ is a
$\fspe(V)$ (and even $\fpe(V)$) invariant and, clearly, it is highest
with respect to $\fb _{+}(W)\oplus \fgl (W)_{-}$.  But such is also
the invariant
$$
\det (v_{i}, v_{\bar j})_{i, j\in I_\ev}
\prod \limits_{s\le t,\; s, t\in I_\ev}(v_{\bar s},
v_{\bar t}).
$$
Set $d=\det (v_{i}, v_{j})_{i, j\in I_\ev}$ and $\Pi ^{-}= \prod
\limits_{s\le t,\; s, t\in I_\ev}(v_{\bar s}, v_{\bar t})$.  This
implies that $X^{+}n=cd\Pi ^{-}$ for $c\neq 0$.  Now, consider the
vector $\omega ^{*k}n$.  By similar arguments the expression
$$
X^{+}(\omega ^{*k}n)=\omega ^{*k}X^{+}n=\omega ^{*k}d\Pi ^{-}c
$$
is an $\fspe(V)$-invariant.  Dividing by $d$ we get the statement of
Lemma.
\end{proof}

Lemma 7.1  implies that on the group $\Pe(V\otimes A)$ there exists a
multiplicative function
$$
B : B^{2}(g)=\Ber~ (g)\; \text { for }\; g\in \Pe(V\otimes A).
$$
Let us denote this function by $\sqrt{\Ber}(g)$ and denote by
$\SPe(V\otimes A)$ the subgroup of $\Pe(V\otimes A)$ consisting of
matrices $g$ such that $\Ber (g)=1$ and denote by $\SSPe(V\otimes A)$
the subgroup consisting of $g$ such that $\sqrt{\Ber}(g)=1$.

$\SSPe(V\otimes A)$ is a connected Lie group corresponding to the Lie
superalgebra $\fspe(V)$.

\ssec{7.3.  Proof of Theorem 1.5} Let us find $\lambda $ for which
there exists a $\fspe(V)$-invariant of type $\lambda $ (we consider
invariants distinct from inner products).  Let $\lambda $ be atypical.
Then an invariant of type $\lambda $, if any, belongs to the algebra
$S^{\bcdot}(V^{*n-1}\oplus \Pi (V)^{*n-1})$ and there exists $k\in \Zee$ such
that
$$
f(g\fL)=(\sqrt{\Ber}(g))^{k}f(\fL).
$$
For generic vectors there exists $g\in \Pe(V\otimes A)$ such that
$$
gv_{i}\in \Span(e_{1}, \dots , e_{n-1}, e_{\od}, \dots ,
e_{\overline{n-1}}).
$$
Then
$$
\Span(g^{-1}e_{n}, g^{-1}e_{\bar n})\perp \Span(e_{1}, \dots , e_{n-1},
e_{\od}, \dots , e_{\overline{n-1}}).
$$
By applying an appropriate transformation from $\Pe(1)$:
$$
\langle g^{-1}e_{n}, g^{-1}e_{\bar n}\rangle\tto\langle e_{n},
e_{\bar n}\rangle
$$
we may assume that $g\in \SSPe(V\otimes A)$.

\centerline{\epsfbox{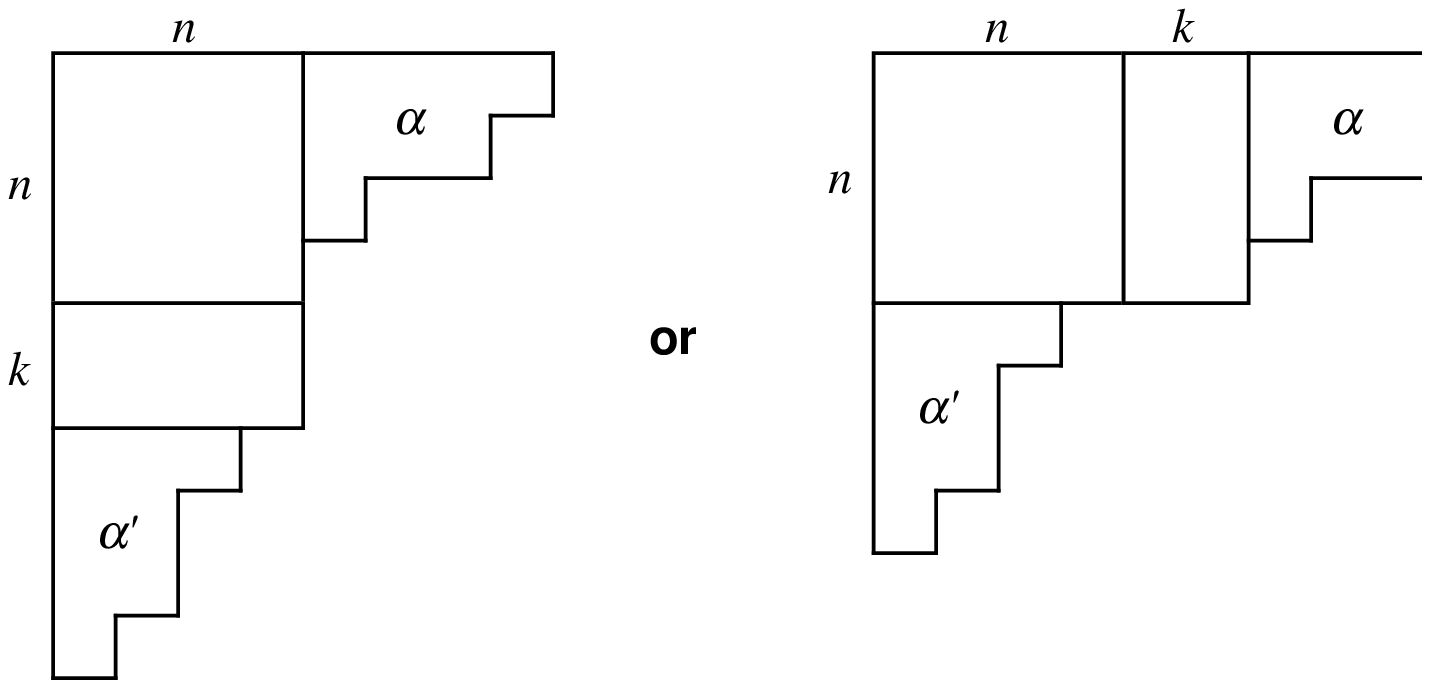}}
\centerline{Figure 7.3}

Let $he_{n}=ae_{n}$ and $he_{\bar n}=a^{-1}e_{\bar n}$, the
other vectors being fixed; then
$$
h\in \Pe(V\otimes A)\; \text{ and }\; \sqrt{\Ber}(h)=a.
$$
Besides, $f(hg\fL)=a^{k}f(g\fL)=f(g\fL)$ and therefore if $k\neq 0$
then $f=0$.  Thus, $\lambda $ should be typical.  By applying Lemma
6.1 we deduce that there exists $k>0$ and $\lambda $ is of the form
shown on Fig.  7.3.

Finally, let us construct invariants of type $\lambda $.  Let $\varphi
_{\alpha }$ be as in Theorem 6.1.  Then $\omega ^{*k}\Pi ^{-}\varphi
_{\alpha }$ and $\Delta ^{*k}\Pi ^{+}\varphi _{\alpha }$ are the
desired invariants.  Theorem is proved.
\qed

\section*{\S 8. Invariants of the Lie superalgebra $\fq (V)$}

Denote by $\GQ(V\otimes A)$ the subgroup in $\GL(V\otimes A)$ that
preserves the inner product
$$
[v^{*}_{t}, v_{s}]=\sum^{n}_{i=1}(x_{ti}x^{*}_{\bar is}+
x_{t\bar i}x^{*}_{is}) \; \text{ for }\; s, t\in I_\ev.
$$
As was mentioned in Theorem 2.2, it suffices to
confine to the algebra $\fA^{n}_{n}=S^{\bcdot}(V^{n}\oplus V^{*n})$.

\ssec{8.1. Proof of Theorem 1.6} Let us consider a
generic collection
$$
\fL=(v_{1}, \dots , v_{n}, v^{*}_{\od}, \dots ,  v^{*}_{\bar n}).
$$
There exists $g\in \GQ(V\otimes A)$ such that $gv_{i}=l_{i}, \; i\in
I_\ev$. If $f$ is a $\GQ(V\otimes A)$-invariant, then
$$
f(\fL)=f(g\fL)=f(e_{1}, \dots , e_{n}, gv^{*}_{1}, \dots , gv^{*}_{n})
$$
is a polynomial in coordinates of $gv^{*}_{t}$ but $(gv^{*}_{t},
e_{s})=(v^{*}_{t}, v_{s})$ and $[gv^{*}_{t}, e_{s}]= [v^{*}_{t},
v_{s}]$; hence, $f$ is a polynomial in inner products.

\section*{\S 9. The invariants of the Lie superalgebra $\fsq (V)$}

First, let us prove a theorem which for the Lie superalgebra $\fq (V)$
plays the same role as Theorem 2.1 plays for $\fgl (V)$.

\ssbegin{9.1}{Theorem} Let $\dim V=(n, n)$, $\dim U=(l, l)$. Then we
have an isomorphism of $\fq (U)\oplus \fq (V)$-modules:
$$
S^{\bcdot}(2^{-1}U\otimes V)\simeq\mathop{\oplus}\limits _{\lambda
_{n+1}=0} 2^{-\delta (| \lambda | )} U^{\lambda }\otimes V^{\lambda }.
$$
Here $U^{\lambda }$ and $V^{\lambda }$ are irreducible $\fq (U)$- and
$\fq (V)$-modules corresponding to $\lambda $ and $\lambda $ is a
strict partition such that $\delta (|\lambda|)= \left\{\begin{matrix}
0&\text{if $|\lambda|$ is even}\cr 1 &
\text{otherwise.}\end{matrix}\right .$ {\em (For the definition of the
module $2^{-1}U\otimes V$ see \cite{S1}.)}
\end{Theorem}

\begin{proof} According to \cite{S1} we have
$$
V^{*\otimes k}=\mathop{\oplus}\limits_{\lambda:
\lambda_{n+1}=0}V^{*\lambda}\otimes T^{\lambda }\cdot 2^{-\delta
(|\lambda|)}, \quad U^{\otimes k}= \mathop{\oplus}\limits_{\mu:
\mu_{n+1}=0}2^{-\delta|\mu|} U^{\mu }\otimes T^{\mu }
$$
Hence,
\begin{align*}
&S^{k}(2^{-1}U\otimes V)=S^{k}(2^{-1}U\otimes (V^{*})^{*})=
S^{k}(\Hom_{G_{1}}(V^{*}, U))= \cr
&\Hom _{G_{k}}(V^{*\otimes k},
U^{\otimes k})=\cr &\mathop{\oplus}\limits _{\lambda , \mu }2^{-\delta
(|\lambda|)}2^{-\delta (| \mu | )}\Hom (V^{*\lambda }, U^{\mu})\otimes
\Hom _{G_{k}}(T^{\lambda }, T^{\mu })=\cr
&\mathop{\oplus}\limits _{\lambda _{n+1}=0}2^{-\delta
(|\lambda|)}U^{\lambda }\otimes V^{\lambda },\text{ here
$G_{k}=\fS_{k}\circ C_{k}$, see \cite{S1}}.
\end{align*}
Theorem is proved.
\end{proof}

\begin{Corollary} We have an isomorphism of
$\fq (V)-\fq (W)-\fq (U)$ trimodules:
$$
S^{\bcdot}(2^{-1}U\otimes V+2^{-1}V^{*}\otimes W)\simeq
\mathop{\oplus}\limits_{\lambda , \mu }2^{-\delta
(|\lambda|)}U^{\lambda } \otimes V^{\lambda }\otimes 2^{-\delta
(|\mu|)} V^{*\mu }\otimes W^{\mu }\simeq\fA^{n}_{n}.
$$
\end{Corollary}

Set
$$
(\fA^{n}_{n})_{\lambda , \mu }=2^{-\delta
(|\lambda|)}U^{\lambda}\otimes V^{\lambda }\otimes 2^{-\delta
(|\mu|)}V^{*\mu }\otimes W^{\mu}
$$
and refer to the elements of this module as the {\it elements of type}
$(\lambda , \mu )$.  The invariants of type $(\lambda , \mu )$ will be
called {\it typical} ones if $\lambda _{n}>0$.

\ssec{9.2} The validity of the following lemma is not difficult to
establish.

\begin{Lemma} Let $\fg =\fq (V)$, $\fsq (V)$ and $\fh $ be a Cartan
subalgebra in $\fg $, let $\fg _{+}$ be the linear span of positive
roots and $L$ the finite dimensional $\fg $-module generated by
$L^{\fg _{+}}$.

Then $L$ is an irreducible $\fg $-module if and only if $L^{\fg _{+}}$
is irreducible as $\fh $-module.
\end{Lemma}

\ssbegin{9.3}{Lemma} Let $\lambda $ and $\mu $ be strict partitions
and $\lambda _{n+1}=\mu _{n+1}=0$.  Then:

{\em a)} $\dim (V^{\lambda }\otimes V^{*\mu })^{\fsq (V)}=0$
if $\lambda \neq \mu$

{\em b)} $\dim(2^{-\delta (|\lambda|)}V^{\lambda }\otimes V^{*\lambda
})^{ \fsq (V)}=1$ if $\lambda _{n}=0$

{\em c)} $\dim (2^{-\delta (|\lambda|)}V^{\lambda }\otimes
V^{*\lambda })=2$ if $\lambda _{n}>0$.
\end{Lemma}

\begin{proof}
$$
(V^{\lambda }\otimes V^{*\mu })^{\fsq (V)}=\Hom _{\fsq (V)}
(V^{\mu }, V^{\lambda }).
$$
Since the even parts of Cartan subalgebras of $\fq (V)$ and $\fsq (V)$
are the same, then the modules $V^{\mu }$ and $V^{\lambda }$ are
non-isomorphic as $\fsq (V)$-modules for $\lambda \neq \mu $ and
heading a) is proved.

Let $\lambda =\mu $ and $\lambda _{n}=0$.  Then by Lemma 9.1
$(V^{\lambda })^{\fq _{+}(V)}$ is an irreducible module and by
\cite{S3} it is of the form
$$
(V^{\lambda })^{\fq _{+}(V)}=\ind^{\fh }_{P_{\lambda }}(k),
$$
where $P_{\lambda }$ is the polarization subordinate to the
functional $\lambda$.

For $\fh =\Span(e_{1}, \dots , e_{n}, e_{\od}, \dots , e_{\bar
n})$ set $\fsh =\Span(e_{1}, \dots , e_{n}, e_{\bar
1}-e_{\bar n}, \dots , e_{\overline{n-1}}-e_{\bar n})$.

Since $\lambda _{n}=0$, then $e_{\bar n}$ belongs to
the kernel of the form
$$
b_{\lambda }(f_{1}, f_{2})=\lambda ([f_{1}, f_{2}]),
\text { where } f_{1}, f_{2}\in \fh _{\od}.
$$
Therefore, the restriction of the form $b_{\lambda }$ to
$(\fsh)_{\od}$ is of the same rank as $b_{\lambda }$.  Therefore, the
module $(V^{\lambda })^{\fq _{+}(V)}$ remains irreducible as
$\fsh$-module and the type of its irreducibility (G or Q) is the same
as that of the $\fh$-module.  This proves heading b).

Let $n$ be even and $\lambda _{n}>0$.  Set $f=\sum\limits^{n}_{i=1}
\frac{e_{\bar i}}{\lambda _{i}}$.  Then we can verify that $b_{\lambda
}(f, f)\neq 0$ and $f\perp (\fsh )_{\od}$.  This proves that the
restriction of the form $b_{\lambda }$ to $(\fsh )_{\od}$ is invariant
but since $\dim (\fsh )_{\bar 1}=n-1$ is odd then, as the
$\fsh$-module, $(V^{\lambda })^{\fq _{+}(V)}$ is irreducible of type
$Q$.

If $n$ is odd, then $(V^{\lambda })^{\fq _{+}(V)}$ is of type $Q$ as
the $\fh$-module.

Since the restriction of $b_{\lambda }$ to $(\fsh )_{\od}$ is
nondegenerate and of an even rank, then $(V^{\lambda })^{\fq _{+}(V)}$
is the direct sum $I\oplus \pi I$, where $I$ is an irreducible
$\fsh$-module of type $G$.

In other words, for $\lambda _{n}>0$ in the module $2^{-\delta
(|\lambda|)}V^{\lambda }\otimes V^{*\lambda}$ there exists an
additional $\fsq (V)$-invariant.  This proves heading b) and Lemma.
\end{proof}

\ssbegin{9.4}{Lemma} Let $\fsq (V)=\fsq (V)\otimes \Span(F)$.  If
$\varphi $ is a typical $\fq (V)$-invariant, then there exists a
unique typical $\fsq (V)$-invariant $\psi $ such that $\varphi =F\psi$.
\end{Lemma}

\begin{proof} By Lemma 9.3 for $\lambda _{n}>0$ there are two
invariants in the module $2^{-\delta (|\lambda|)}V^{\lambda } \otimes
V^{*\lambda }$: one is an $\fq (V)$-invariant $\varphi $ and another
one is a $\fsq (V)$-invariant $\psi $.  Hence, $F\psi \neq 0$ is,
clearly, a $\fq (V)$-invariant; hence, $F\psi =c\varphi , \; c\in
\Cee$.  Setting $\psi =\psi /c$ we get the statement of Lemma.
\end{proof}

\ssbegin{9.5}{Lemma} Any $\fsq (V)$-invariant which is not a $\fq
(V)$-invariant is of the form $\varphi\qet Y$, where $\varphi$ is a
$\fq (V)$-invariant and $Y$ is the same as in Theorem 1.7.
\end{Lemma}

\begin{proof} Let us take a Grassmann algebra $A$ with sufficiently
large number of generators and consider the elements of the algebra
$\fA^{n}_{n}$ as functions on the space of collections
$$
\fL=(v_{1}, \dots
, v_{n}, v^{*}_{1}, \dots
, v^{*}_{n}), \text {where}\; v_{i}\in (V\otimes A)_\ev\;
\text {and}\;
v^{*}_{i}\in (\Hom _{A}(V\otimes A, A))_\ev.
$$
Let $f$ be an $\fsq (V)$-invariant which is not a $\fq (V)$-invariant
and let $M$ be the set of collections $\fL$ such that $\{v_{1}, \dots
, v_{n}\}$ is a basis in $(V\otimes A)_\ev$.  Denote by $SQ(V\otimes
A)$ the subgroup of transformations from $\GQ(V\otimes A)$ whose queer
determinant is equal to 1.  Take $g\in \GQ(V\otimes A)$ such that
$$
ge_{i}=v_{i}\; \text {and $he_{i}=e_{i}+e_{\bar i}\xi , he_{\bar
i}=e_{\bar i}+e_{i}\xi $, where $\xi =\qet\; g/n$.}
$$
Then $hg^{-1}\in SQ(V\otimes A)$ and
$$
f(\fL)=f(hg^{-1}\fL)=f(he_{1}, \dots , he_{n}, hg^{-1}v^{*}_{1}, \dots ,
hg^{-1}v^{*}_{n})
$$
is a polynomial in $\xi $ and coordinates of $hg^{-1}v^{*}_{i}$.  But
$\xi =\qet g/n=\qet Y/n$ and
$$
(hg^{-1}v^{*}_{i}, e_{j})=(v^{*}_{i}, gh^{-1}e_{j})=(v^{*}_{i},
ge_{j}-ge_{\bar j}\xi )=(v^{*}_{i}, v_{j})-\langle v^{*}_{i},
v_{j}\rangle \xi .
$$
Lemma is proved.
\end{proof}

\ssbegin{9.6}{Lemma} Let $\varphi $ be a $\fq (V)$-invariant.  Then
$\varphi \qet Y$ is a polynomial if and only if $\varphi $ is a
typical invariant.
\end{Lemma}

\begin{proof} First, let us prove that (in notations of Lemma 9.4) if
$\qtr F=1$, then $F (\qet Y)=1$.

Indeed, let $h$ be selected as in Lemma 9.5. Then
\begin{align*}
&\qet Y+F\xi \qet Y=\qet (\exp (F\xi )Y)= \cr &=\qet (\exp (F\xi
))+\qet Y=\xi +\qet Y;
\end{align*}
hence, $F(\qet Y)=1$.  Let $\varphi $ be a typical $\fq
(V)$-invariant.  Then by Lemma 9.4 there exists a unique $\fsq
(V)$-invariant $\psi $ such that $\varphi =F\psi $.  On the other
hand, by Lemma 9.5 $\psi =\varphi _{1}\qet Y$, where $\varphi _{1}$ is
a $\fq (V)$-invariant.  Hence,
$$
\varphi =F\psi =F(\varphi _{1}\qet Y)= \pm \varphi _{1}F(\qet Y)=\pm
\varphi _{1}
$$
and, therefore, $\varphi \qet Y= \pm \psi _{1}$ is a polynomial.
Since it is $\fsq (V)$-invariant, then by Lemma 9.2 it is a typical
one and therefore $\varphi =F(\varphi \qet Y)$ is also a typical
invariant.  Lemma is proved.
\end{proof}

The above arguments show that to construct $\fsq (V)$-invariants it
suffices to construct typical $\fq (V)$-invariants.  One of the ways
of such a construction is described in the following lemma.

\ssbegin{9.7}{Lemma} {\em (Notations from Theorem 1.7.)} For any
partition $\lambda $ such that $\lambda _{1}>\dots >\lambda _{n}>0$
the polynomial
$$
p_{\lambda }=\qtr Z^{\lambda _{1}}\dots
\qtr Z^{\lambda _{n}}
$$
is a typical $\fq (V)$-invariant.
\end{Lemma}

{\bf An idea of the proof.}
$$
(\fA^{n}_{n})^{\fq (V)}=S^{\bcdot}(2^{-1}(U\otimes W))
$$
is a $\fq (U)\oplus \fq (W)$-module such that $\dim U=\dim W=(n, n)$.
Take a superspace $L$ such that $\dim L=(n, n)$ and fix isomorphisms
$L=U$ and $L^{*}=\pi (W)$ which determine isomorphisms of algebras
$S^{\bcdot}(\fq (L)^{*})=(\fA^{n}_{n})^{\fq (V)}$.  For an irreducible
representation $\pi$ the functionals $\str \pi $ (and $\qtr \pi$ if
the representation is of type $Q)$ are $\fq (V)$-invariant elements of
the algebra $S^{\bcdot}(\fq (L))^{*}$ and if $\pi _{\lambda }$ corresponds to
the irreducible module $U^{\lambda }$ then $\qtr \pi _{\lambda }$ (or
$\str \pi _{\lambda })$ restricted to $S^{| \lambda | }(\fq (L)^{*})$
is of type $\lambda $.  The invariant elements are uniquely determined
by their restrictions to a Cartan subalgebra in $\fq (L)$.

It is not difficult to verify that $\str \pi _{\lambda }$ (or $\qtr
\pi _{\lambda })$ and $p_{\lambda }$ have the identical restrictions.

\ssec{9.8.  Proof of Theorem 1.7} Lemma 9.7 provides with a
construction of a typical $\fq (V)$-invariant and Lemma 9.6 with the
construction of an $\fsq (V)$-invariant of type $\lambda$ which
completes the proof of the theorem.

\section*{Appendix 0. Background}

\ssec{A0.1.  Linear algebra in superspaces.  Generalities}

A {\it superspace} is a $\Zee /2$-graded space; for a superspace
$V=V_{\ev}\oplus V_{\od }$ denote by $\Pi (V)$ another copy of the
same superspace: with the shifted parity, i.e., $(\Pi(V))_{\bar i}=
V_{\bar i+\od }$.  The {\it superdimension} of $V$ is $\dim
V=p+q\varepsilon$, where $\varepsilon^2=1$ and $p=\dim V_{\ev}$,
$q=\dim V_{\od }$.  (Usually $\dim V$ is expressed as a pair $(p, q)$
or $p|q$; this obscures the fact that $\dim V\otimes W=\dim V\cdot
\dim W$ which is clear with the use of $\varepsilon$.)

A superspace structure in $V$ induces the superspace structure in the
space $\End (V)$.  A {\it basis of a superspace} always
consists of {\it homogeneous} vectors; let $\Par=(p_1, \dots,
p_{\dim V})$ be an ordered collection of their parities.  We call
$\Par$ the {\it format} of the basis of $V$.  A square {\it
supermatrix} of format (size) $\Par$ is a $\dim V\times \dim V$ matrix
whose $i$th row and $i$th column are of the same parity $p_i$.  The
matrix unit $E_{ij}$ is supposed to be of parity $p_i+p_j$ and the
bracket of supermatrices (of the same format) is defined via Sign
Rule:

{\it if something of parity $p$ moves past something of parity
$q$ the sign $(-1)^{pq}$ accrues; the formulas defined on homogeneous
elements are extended to arbitrary ones via linearity}.

Examples of application of Sign Rule:
setting $[X, Y]=XY-(-1)^{p(X)p(Y)}YX$ we get the notion of the
supercommutator and the ensuing notions of the supercommutative
superalgebra and the Lie superalgebra (that in addition to
superskew-commutativity satisfies the super Jacobi identity, i.e., the
Jacobi identity amended with the Sign Rule). The superderivation of a
superalgebra $A$ is a linear map $D: A\tto A$ such that satisfies the
Leibniz rule (and Sign rule)
$$
D(ab)=D(a)b+(-1)^{p(D)p(a)}aD(b).
$$

Usually, $\Par$ is of the form $(\ev , \dots, \ev , \od , \dots,
\od)$.  Such a format is called {\it standard}.  In this paper we can
do without nonstandard formats.  But they are vital in the study of
systems of simple roots that the reader might be interested in.

The {\it general linear} Lie superalgebra of all supermatrices of size
$\Par$ is denoted by $\fgl(\Par)$; usually, $\fgl(\ev, \dots, \ev,
\od, \dots, \od)$ is abbreviated to $\fgl(\dim V_{\bar 0}|\dim V_{\bar
1})$.  Any matrix from $\fgl(\Par)$ can be expressed as the sum of its
even and odd parts; in the standard format this is the block
expression:
$$
\begin{pmatrix}A&B\\ C&D\end{pmatrix}=\begin{pmatrix}A&0\\
0&D\end{pmatrix}+\begin{pmatrix}0&B\\ C&0\end{pmatrix},\quad
p\left(\begin{pmatrix}A&0\\
0&D\end{pmatrix}\right)=\ev, \; p\left(\begin{pmatrix}0&B\\
C&0\end{pmatrix}\right)=\od.
$$

The {\it supertrace} is the map $\fgl (\Par)\longrightarrow \Cee$,
$(A_{ij})\mapsto \sum (-1)^{p_{i}}A_{ii}$.  Since $\str [x, y]=0$, the
space of supertraceless matrices constitutes the {\it special linear}
Lie subsuperalgebra $\fsl(\Par)$.

There are, however, two super versions of $\fgl(n)$, not one.  The
other version is called the {\it qeer} Lie superalgebra and is defined
as the one that preserves the complex structure given by an {\it odd}
operator $J$, i.e., is the centralizer $C(J)$ of $J$:
$$
\fq(n)=C(J)=\{X\in\fgl(n|n): [X, J]=0 \}, \text{ where } J^2=-\id.
$$
It is clear that by a change of basis we can reduce $J$ to the form
$J_{2n}=\begin{pmatrix}0&1_n\\ -1&0\end{pmatrix}$.  In the standard
format we have
$$
\fq(n)=\left \{\begin{pmatrix}A&B\\ B&A\end{pmatrix}\right\}.
$$
On $\fq(n)$, the {\it queer trace} is defined: $\qtr:
\begin{pmatrix}A&B\\
B&A\end{pmatrix}\mapsto
\tr B$. Denote by $\fsq(n)$ the Lie superalgebra of {\it queertraceless}
matrices.

Observe that the identity representations of $\fq$ and $\fsq$ in $V$,
though irreducible in supersence, are not irreducible in the
nongraded sence: take homogeneous linearly independent vectors $v_1$,
\dots , $v_n$ from $V$; then $\Span (v_1+J(v_1), \dots , v_n+J(v_n))$
is an invariant subspace of $V$ which is not a subsuperspace.

A representation is {\it irreducible} \index{representation of Lie
superalgebra irreducible} \index{$G$-type irreducible representation
of Lie superalgebra} of {\it general type} or just of {\it $G$-type}
if there is no invariant subspace, otherwise it is called {\it
irreducible of $Q$-type} \index{$Q$-type irreducible representation of
Lie superalgebra} ($Q$ is after the general queer Lie superalgebra ---
a specifically superish analog of $\fgl$); an irreducible
representation of $Q$-type has no invariant sub{\it super}space but
{\it has} an invariant subspace.

{\bf Superalgebras that preserve bilinear forms: two types}.  To the
linear map $F$ of superspaces there corresponds the dual map $F^*$
between the dual superspaces; if $A$ is the supermatrix corresponding
to $F$ in a basis of the format $\Par$, then to $F^*$ the {\it
supertransposed} matrix $A^{st}$ corresponds:
$$
(A^{st})_{ij}=(-1)^{(p_{i}+p_{j})(p_{i}+p(A))}A_{ji}.
$$

The supermatrices $X\in\fgl(\Par)$ such that
$$
X^{st}B+(-1)^{p(X)p(B)}BX=0\quad \text{for an homogeneous matrix
$B\in\fgl(\Par)$}
$$
constitute the Lie superalgebra $\faut (B)$ that preserves the
bilinear form on $V$ with matrix $B$.  Most popular is the
nondegenerate supersymmetric form whose matrix in the standard format
is the canonical form $B_{ev}$ or $B'_{ev}$:
$$
B_{ev}(m|2n)= \begin{pmatrix}
1_m&0\\
0&J_{2n}
\end{pmatrix},\quad \text{where
$J_{2n}=\begin{pmatrix}0&1_n\\-1_n&0\end{pmatrix}$,}
$$
or
$$
B'_{ev}(m|2n)= \begin{pmatrix}
\antidiag (1, \dots , 1)&0\\
0&J_{2n}
\end{pmatrix}.
$$
The usual notation for $\faut (B_{ev}(m|2n))$ is $\fosp(m|2n)$ or
$\fosp^{sy}(m|2n)$.  (Observe that the passage from $V$ to $\Pi (V)$
sends the supersymmetric forms to superskew-symmetric ones, preserved
by the \lq\lq symplectico-orthogonal" Lie superalgebra $\fsp'\fo
(2n|m)$ or $\fosp^{sk}(m|2n)$ which is isomorphic to
$\fosp^{sy}(m|2n)$ but has a different matrix realization.  We never
use notation $\fsp'\fo (2n|m)$ in order not to confuse with the
special Poisson superalgebra.)

In the standard format the matrix realizations of these algebras
are:
$$
\begin{matrix}
\fosp (m|2n)=\left\{\left (\begin{matrix} E&Y&X^t\\
X&A&B\\
-Y^t&C&-A^t\end{matrix} \right)\right\};\quad \fosp^{sk}(m|2n)=
\left\{\left(\begin{matrix} A&B&X\\
C&-A^t&Y^t\\
Y&-X^t&E\end{matrix} \right)\right\}, \\
\text{where}\;
\left(\begin{matrix} A&B\\
C&-A^t\end{matrix} \right)\in \fsp(2n),\quad E\in\fo(m)\;
\text{and}\;  {}^t \; \text{is the usual transposition}.\end{matrix}
$$

A nondegenerate supersymmetric odd bilinear form $B_{odd}(n|n)$ can be
reduced to the canonical form whose matrix in the standard format is
$J_{2n}$.  A canonical form of the superskew odd nondegenerate form in
the standard format is $\Pi_{2n}=\begin{pmatrix}
0&1_n\\1_n&0\end{pmatrix}$.  The usual notation for $\faut
(B_{odd}(\Par))$ is $\fpe(\Par)$.  The passage from $V$ to $\Pi (V)$
sends the supersymmetric forms to superskew-symmetric ones and
establishes an isomorphism $\fpe^{sy}(\Par)\cong\fpe^{sk}(\Par)$.
This Lie superalgebra is called, as A.~Weil suggested, {\it
periplectic}.  In the standard format these superalgebras are
shorthanded as in the following formula, where their matrix
realizations is also given:
$$
\begin{matrix}
\fpe ^{sy} (n)=\left\{\begin{pmatrix} A&B\\
C&-A^t\end{pmatrix}, \; \text{where}\; B=-B^t,
C=C^t\right\};\\
\fpe^{sk}(n)=\left\{\begin{pmatrix}A&B\\ C&-A^t\end{pmatrix}, \;
\text{where}\; B=B^t, C=-C^t\right\}.
\end{matrix}
$$

The {\it special periplectic} superalgebra is $\fspe(n)=\{X\in\fpe(n): \str
X=0\}$.

Observe that though the Lie superalgebras $\fosp^{sy} (m|2n)$ and
$\fpe ^{sk} (2n|m)$, as well as $\fpe ^{sy} (n)$ and $\fpe ^{sk} (n)$,
are isomorphic, the difference between them is sometimes crucial.

\ssec{A0.1.1.  Projectivization} If $\fs$ is a Lie algebra of scalar
matrices, and $\fg\subset \fgl (n|n)$ is a Lie subsuperalgebra
containing $\fs$, then the {\it projective} Lie superalgebra of type
$\fg$ is $\fpg= \fg/\fs$.

Projectivization sometimes leads to new Lie superalgebras, for
example: $\fpgl (n|n)$, $\fpsl (n|n)$, $\fpq (n)$, $\fpsq (n)$;
whereas $\fpgl (p|q)\cong \fsl (p|q)$ if $p\neq q$.

\section*{Appendix 1.  Certain constructions with the point functor} The
point functor is well-known in algebraic geometry since at least 1953
\cite{Wi}.  The advertising of ringed spaces with nilpotents in the
structure sheaf that followed the discovery of supersymmetries caused
many mathematicians and physicists to realize the usefulness of the
language of points.  F.~A.~Berezin \cite{B} was the first who applied the
point
functor to study Lie superalgebras.  Here we present some of his
results and their generalizations.

All superalgebras and modules are supposed to be
finite dimensional over $\Cee $.

\ssec{A0.  What is a Lie superalgebra} A {\it Lie
superalgebra}\index{Lie superalgebra} $L= L_{\ev} \oplus L_{\od}$ is a
linear supermanifold $\cL=(\cL_{rd}, \cO _{\cL})$ such that for
\lq\lq any" (say, finitely generated, or from some other appropriate
category) supercommutative superalgerba $C$ the space $\cL(C)=\Hom
(\Spec C, \cL)$, called {\it the space of $C$-points of} $\cL$ is a
Lie algebra and the correspondence $C\longrightarrow \cL(C)$ is a
functor in $C$.

\begin{rem*}{Exercise} What is the set of $C$-points of a Lie
superalgebra $L$ given in the conventional definition as a superspace
with a superskew-symmetric product satisfying super Jacobi identity?
What is the supermanifold $\cL=(\cL_{rd}, \cO _{\cL})$ that
corresponds to $L$?  (Answer: $(L \otimes C)_{\ev }$ and
$\cL_{rd}=L_{\ev}$, $\cO_{\cL}=\cO _{L_{\ev}}\otimes
S^{\bcdot}(L_{\od}^*)$.  Recall that $S^{\bcdot}$ applied to the odd
superspace is understood in the supersence, i.e., as the ordinary
$\Lambda^{\bcdot}$ applied to this space.)
\end{rem*}

A {\it Lie superalgebra homomorphism} $\rho: L_1 \tto L_2$ in these
terms is a functor morphism, i.e., a collection of Lie algebra
homomorphisms $\rho_C: L_1 (C)\tto L_2(C)$ compatible with morphisms
of supercommutative superalgebras $C\longrightarrow C'$.

In particular, a {\it representation}\index{representation of Lie
superalgebra} of a Lie superalgebra $L$ in a superspace $V$ is a
homomorphism $\rho: L\longrightarrow \fgl (V)$, i.e., a collection of
Lie algebra homomorphisms $\rho_C: L(C) \longrightarrow ( \fgl (V
\otimes C))_{{\ev }}$.

Example. $\qtr$ is not a representation of $\fq(n)$ according to the naive
definition (``a representation is a Lie superlagebra homomorphism'',
hence, an even map), but is a representation, moreover, an irreducible
one, if we consider odd parameters.

Thus, let $\fg$ be a Lie superalgebra, $V$ a $\fg$-module, $\Lambda $
the Grassmann superalgebra over $\Cee $ generated by $q$
indeterminates.  Define $\varphi :\Lambda \otimes V^{*}\tto \Hom
_{\Lambda }(\Lambda \otimes V, \Lambda )$ by setting
$$
\varphi (\xi \otimes \alpha )(\eta \otimes v)=(-1)^{p(\alpha )(\eta
)}\xi \eta \alpha (v), \; \text{for any $\xi, \eta \in \Lambda , \alpha
\in V^{*}$.}
$$
Extend the ground field to $\Lambda $ and consider $\Lambda \otimes
V^{*}$ and $\Hom_\Lambda(\Lambda \otimes V, \Lambda )$ as $\Lambda
\otimes \fg$-modules.

\ssbegin{A1}{Lemma} $\varphi $ is a $\Lambda \otimes \fg$-module
isomorphism.
\end{Lemma}

{\it Proof}.  Since $V$ is finite dimensional, $\varphi $ is a vector
space isomorphism over $\Lambda $; besides, it is obvious that
$\varphi $ is a $\Lambda $-module homomorphism.  Now take
$$
\xi _{1}, \xi _{2}, \xi _{3}\in \Lambda,\; \; \alpha \in V^{*},\;
\; v\in V, \; \; x\in \fg.
$$
It is an easy exercise to prove that
$$
[(\xi _{1}\otimes x)\varphi (\xi _{2}\otimes \alpha )](\xi _{3}\otimes
v)=\varphi (\xi _{1}\otimes x(\xi _{2}\otimes \alpha ))(\xi _{3}\otimes
v).  \qed
$$

Consider the composition of maps
$$
V^{*}\buildrel{\varphi _{1}}\over\tto\Lambda \otimes
V^{*}\buildrel{\varphi }\over\tto\Hom _{\Lambda }(\Lambda
\otimes V, \Lambda )\buildrel{\varphi _{2}}\over\tto S_{\Lambda
}(\Hom _{\Lambda }(\Lambda \otimes V, \Lambda )),
$$
where $\varphi _{1}(\alpha )=1\otimes \alpha$ and $\varphi _{2}$ is a
canonical embedding of a module in its symmetric algebra.  The $\Cee
$-module homomorphism $\varphi _{2}\circ \varphi \circ \varphi _{1}$
induces the algebra homomorphism
$$
S(V^{*})=S_{\Cee }(V^{*})\tto S_{\Lambda }(\Hom _{\Lambda
}(\Lambda \otimes V, \Lambda ))
$$
and, since the latter algebra is a $\Lambda $-module, we get an
algebra homomorphism
$$
\Lambda \otimes S(V^{*})\buildrel{\psi }\over\tto S_{\Lambda }(\Hom
_{\Lambda }(\Lambda \otimes V, \Lambda )).
$$
Besides, both algebras possess a natural $\Lambda \otimes
\fg$-module structure.

\ssbegin{A2}{Lemma}  $\psi $ {\it is a $\Lambda \otimes
\fg$-modules and $\Lambda \otimes \fg$-algebras isomorphism}.
\end{Lemma}

\begin{proof} Let us construct the inverse homomorphism.
Consider the composition
$$
\Hom _{\Lambda }(\Lambda \otimes V, \Lambda )\buildrel{\varphi
^{-1}}\over\tto \Lambda \otimes V^{*}\tto\Lambda
\otimes S(V^{*}).
$$
Since this composition is a $\Lambda $-module homomorphism, it
induces
the homomorphism
$$
\tilde{\psi}: S_{\Lambda }(\Hom _{\Lambda } (\Lambda \otimes V,
\Lambda ))\tto\Lambda \otimes S(V^{*}).
$$
It is not difficult to verify that
$$
\psi \circ \tilde{\psi }|_{\Hom _{\Lambda }(\Lambda \otimes V, \Lambda
)}=\id;\qquad \tilde{\psi }\circ \psi |_{\Lambda \otimes
S(V^{*})}=\id;
$$
hence, $\psi$ is an isomorphism and $\tilde{\psi }$ is its inverse.
The following proposition shows that $\psi $ is a $\Lambda \otimes
\fg$-module isomorphism and completes the proof of Lemma A2.
\end{proof}

\ssbegin{A3}{Proposition} Let $A$, $B$ be $\Lambda $-superalgebras,
$\fg$ a Lie superalgebra over $\Lambda $ acting by differentiations on
$A$ and $B$.  Let $M\subset A$, $N\subset B$ be $\Lambda $-submodules
which are at the same time $\fg$-modules generating $A$ and $B$,
respectively, $f: A\tto B$ an algebra homomorphism such that
$f(M)\subset N$ and $f|_{M}$ is a $\fg$-module homomorphism.  Then $f$
is a $\fg$-module homomorphism.
\end{Proposition}

\begin{proof} Let $a\in A$.  We may assume that $a=a_{1}\dots a_{n}$,
where the $a_{i}\in M$.  Then for $x\in \fg$ we have
$$
\renewcommand{\arraystretch}{1.2}
\begin{gathered}
f(x(a_{1}\dots a_{n}))=f(\sum \pm a_{1}\dots xa_{i}\dots a_{n})=
\sum \pm f(a_{1})\dots f(xa_{i})\dots f(a_{n})\\
=\sum \pm f(a_{1})\dots
xf(a_{i})\dots f(a_{n})= x[f(a_{1})\dots
f(a_{n})]=xf(a_{1}\dots a_{n}).
\end{gathered}
$$

This proves Proposition and completes the proof of
Lemma A2. \end{proof}

Now, let $\fh$ be a Lie superalgebra over $\Lambda $ and $U$ be a
$\Lambda $ and $\fh$-module.  Consider $U_{\ev}$ as a $\Cee $-module.
Then, clearly, the natural embedding $U_{\ev}\tto U$ is extendable to
a $\Lambda $-module homomorphism $\varphi : \Lambda \otimes
U_{\ev}\tto U$.

\ssbegin{A4}{Lemma} The homomorphism $\varphi $ is an
$\fh_{\ev}$-module homomorphism.
\end{Lemma}

\begin{proof} Let $x\in \fh_{\ev}, \xi \in \Lambda $ and $u\in
U_{\ev}$. Then
$$
\renewcommand{\arraystretch}{1.2}
\begin{gathered}
\varphi (x(\xi \otimes u))=\varphi (\xi \otimes xu)=\xi xu, \\
x\varphi (\xi \otimes u)=x\xi u=\text{ ( by definition of a module over a
superalgebra) } \xi xu. \qed
\end{gathered}
$$

Thus, the adjoint map
$$
\Hom_{\Lambda }(U, \Lambda )\tto\Hom_{\Lambda }(\Lambda
\otimes U_{\ev}, \Lambda )
$$
is also an $\fh_{\ev}$-module homomorphism, therefore, by Proposition
A3 the algebra homomorphism
$$
S_{\Lambda}(\Hom _{\Lambda}(U, \Lambda ))\tto
S_{\Lambda}(\Hom_{\Lambda }(\Lambda \otimes U_{\ev}, \Lambda ))
$$
induced by this map is at the same time a $\fh_{\ev}$-module morphism.
Besides, by Lemma A2 the algebra $S_{\Lambda }(\Hom _{\Lambda
}(\Lambda \otimes U_{\ev}, \Lambda))$ is isomorphic as a $\Lambda
\otimes \fh_{\ev}$-module and as an algebra to $\Lambda \otimes
S(U^{*}_{\ev})$.  In particular, they are isomorphic as
$\fh_{\ev}$-modules.

Denote by $\theta$ the composition of the homomorphisms
$$
S(V^{*})\tto\Lambda \otimes S(V^{*})\tto S_{\Lambda}(\Hom_{\Lambda
}(\Lambda \otimes V, \Lambda ))\tto S_{\Lambda}
(\Hom_{\Lambda}(\Lambda \otimes U_{\ev}, \Lambda )),
$$
where $U_{\ev}=(\Lambda \otimes V)_{\ev}=V_{\Lambda }$.
\end{proof}

\ssbegin{A5}{Proposition} If $q>$dimV$_{\od}$ and $\xi \in \Lambda ,
p(\xi )=\od$, then the restriction of $\theta $ onto $\Cee [\xi
]\otimes S(V^{*})$ is injective.
\end{Proposition}

\begin{proof} If $u\in V_{\Lambda}$, then there is defined a linear form
$L_{u}:\Hom _{\Lambda }(\Lambda \otimes V_{\Lambda }, \Lambda
)\tto\Lambda $ by the formulas $L_{u}(l)=l(1\otimes u)$ and
$L_{u}(\xi l)=\xi l(1\otimes u)=\xi L_{u}(l)$.

Therefore, $L_{u}$ is a $\Lambda $-module homomorphism, hence, it is
uniquely extendable to a homomorphism
$$
\varphi _{u}: \fa=S_{\Lambda }(\Hom _{\Lambda }(\Lambda \otimes
V_{\Lambda }, \Lambda ))\tto\Lambda
$$
Consider the elements of $\fa$ as functions on $V_{\Lambda }$ setting
$f(u)=\varphi _{u}(f)$ for $f\in \fa$ and $u\in V_{\Lambda }$.  If
$f\in \Lambda \otimes S(V^{*})$, then set $f(u)=\varphi _{u}\circ
\theta (f)$.  For $\alpha \in V^{*}, \xi \in \Lambda $ we have
$$
(\xi \otimes \alpha )(u)=\varphi _{u}\circ \theta (\xi \otimes \alpha
)=L_{u}\circ \theta (\xi \otimes \alpha )=\theta (\xi \otimes \alpha
)(1\otimes u).
$$
If $\{e_{i}\}_{i\in I}$ is a basis in $V$ and $u=\sum \lambda _{i}\otimes
e_{i}$,
then
$$
(\xi \otimes \alpha )(u)=\sum (-1)^{p(\alpha )p(e_{i})}\xi \lambda
_{i}\alpha (e_{i}).  \eqno{(1)}
$$
On the other hand, the algebra $\Cee [\xi ]\otimes S(V^{*})$ is
identified with the free supercommutative superalgebra generated by
the $e^{*}_{i}$ and $\xi $.

Let
us assume that $p(e^{*}_{i})=0$ for $i\le n$ and $p(e^{*}_{i})=1$ for $i>n$. If
$f\in \Cee [\xi ]\otimes S(V^{*})$, then
$$
f=f_{0}+\xi f_{1}, f_{j}=\sum f_{ji_{1}\dots i_{k}}e^{*}_{i_{1}}\dots
e^{*}_{i_{k}}, \; \text{where}\; j=0, 1\; \text{and}\; f_{ji_{1}\dots
i_{k}}\in S(V^{*}_{\ev}).
$$
By $(1)$ we have
$$
\renewcommand{\arraystretch}{1.2}
\begin{gathered}
f(u)=f_{0}(u)+\xi f_{1}(u)=\\
\sum f_{0i_{1}\dots i_{k}}(u)e^{*}_{i_{1}}(u)\dots
e^{*}_{i_{k}}(u)+\sum f_{1i_{1}\dots i_{k}}(u)e^{*}_{i_{1}}(u)\dots
e^{*}_{i_{k}}(u).
\end{gathered}
$$

Set $\lambda _{i}=a_{i}$ for $i\le n$ and $\lambda _{i}=\xi _{i-n}$
for $i>n$.  Then since $q>\dim V_{\od}$, we may assume that the family
$\{\xi _{i}\}_{i\in I}$ freely generates $S(V^{*}_{\od})$ and
$$
f(u)=\sum (-1)^{k}f_{i_{1}\dots i_{k}}(a_{1}\dots a_{n})\xi
_{i_{1}-n\dots i_{k}-n}. \eqno{(2)}
$$

If $\theta (f)=0$, then $f(u)=\varphi _{u}\circ \theta (f)$ for any
$u\in V_{\Lambda }$.  It follows from $(2)$ that $f_{i_{1}\dots
i_{k}}(a)=0$ for any $a\in \Cee ^{n}$.  But since $\Cee $ is
algebraically closed, it follows (with Prop.  5.3.1 from \cite{Bu})
that $f_{i_{1}\dots i_{k}}=0$; hence, $f=0$.
\end{proof}

\ssbegin{A6}{Lemma} Let $q>$dimV$_{\od}$.  Then $f\in S(V^{*})$ is
a $\fg$-invariant if and only if $\theta (f)\in \Lambda \otimes
S(V^{*}_{\Lambda })$ is $\fg_{\Lambda }$-invariant.
\end{Lemma}

\begin{proof} Consider the factorization of $\theta $:
$$
\begin{gathered}
S(V^{*})\buildrel{i_{1}}\over\tto\Lambda \otimes
S(V^{*})\buildrel{i_{2}}\over\tto\\
S_{\Lambda }(\Hom _{\Lambda
}(\Lambda \otimes V, \Lambda ))\buildrel{i_{3}}\over\tto S_{\Lambda
}(\Hom _{\Lambda }(\Lambda \otimes V_{\Lambda }, \Lambda
))\buildrel{i_{4}}\over\tto\Lambda \otimes S(V^{*}_{\Lambda }).
\end{gathered}
$$
Let $f\in S(V^{*})^\fg$, then
$$
(\xi \otimes x)(i_{1}(f)=(\xi \otimes x)(1\otimes f)=\xi \otimes
xf=0\; \; \text{ for }\; \xi \in \Lambda , x\in \fg.
$$

Conversely, let $yi_{1}(f)=0$ for any $y\in \fg_{\Lambda }$.
Then
$$
0=(\xi \otimes x)(1\otimes f)=\xi \otimes xf.
$$
If $p(y)=\od$ then let $p(\xi )=\od$. Therefore, $\xi \otimes xf=0$ and
$yf=0$.

Thus, conditions
$$
f\in S(V^{*})^\fg\longleftrightarrow i_{1}(f)\in (\Lambda
\otimes S(V^{*}))^\fg_{\Lambda }.
$$
Since $i_{2}$, $i_{3}$, $i_{4}$ are $\fg_{\Lambda }$-module
homomorphisms, the above implies that if $f$ is a $\fg$-invariant then
$\theta (f)=i_{4}\circ i_{3}\circ i_{2}\circ i_{1}(f)$ is also a
$\fg_{\Lambda }$-invariant.

Conversely, let $\theta (f)$ be a $\fg_{\Lambda }$-invariant.  Let
$x\in\fg_{\ev}$.  Then
$$
\theta (1\otimes xf)=\theta ((1\otimes x)(1\otimes f))=(1\otimes
x)\theta (f)=0.
$$

By Proposition A5 $1\otimes xf=0$ and $xf=0$. Let $x\in
\fg_{\od}, \; \xi \in \Lambda , \; p(\xi )=\od$. Then
$\theta (\xi \otimes xf)= (\xi \otimes x)\theta (1\otimes f)=0$ and
again by Proposition A5 $\xi \otimes xf=0$; hence $xf=0$ and
therefore, $f\in S(V^{*})^\fg{{ g}}$.
\end{proof}

\ssec{A7. Remark} The point of the above lemmas and
propositions is that while seeking invariant polynomials on $V$
we may consider them as functions on $V_{\Lambda }$ invariant with
respect to the Lie algebra $\fg_{\Lambda }$. It makes it possible to
apply the theory of usual Lie groups and Lie algebras and their
representations.

\ssec{A8. Remark} Let $\varphi $ be an elementary automorphism
(of the form $\theta _{\beta }$ in Lemma 1.2.3 below) of the Lie
algebra $\fg_{\ev}$. Clearly, $\varphi $ can be uniquely
extended to an automorphism of the Lie superalgebra $\fg$. Let
$\varphi (\fh)=\fh$, where $\fh$ is a Cartan subalgebra of
$\fg$. If $i: S(\fg^{*})^{\fg}\tto S(\fh^{*})$ is the restriction
homomorphism, then, clearly, $i(S(\fg^{*})^\fg)\subset
S(\fh^{*})^{\varphi }$, where $A^{\varphi }$ is the set of $\varphi
$-invariant elements of $A$.

\ssbegin{A9}{Proposition} Let $A$ be a commutative finitely
generated algebra over $\Cee $ without nilpotents, $\fa=A\otimes
\Lambda (p)$. Let $q\ge p$ and $f\in \fa$ be such that
$\varphi (f)=0$ for any $\varphi :\fa\tto\Lambda (q)$. Then
$f=0$.
\end{Proposition}

\begin{proof} Let $\psi: A\tto\Cee $ be an arbitrary homomorphism.
Let us extend $\psi $ to a homomorphism $\varphi : \fa\tto\Lambda (q)$
setting $\varphi =\psi \otimes 1$.  If $\xi _{1}, \dots , \xi _{p}$
are generators of $\Lambda (p), f\in \fa$ and $f=\sum f_{i_{1}\dots
i_{k}}\xi _{i_{1}}\dots \xi _{i_{k}}$, then the condition $\varphi
(f)=0$ yields $\psi (f_{i_{1}\dots i_{k}})=0$ and, since $\psi $ is
arbitrary, then Proposition 5.3.1 in \cite{Bu}  shows that
$f_{i_{1}\dots i_{k}}=0$; hence, $f=0$.
\end{proof}


\begin{thebibliography}{9999}
\bibitem[B]{B}
Berezin F., Representations of the supergroup $U(p, q)$.
Funkcional. Anal. i Prilozhen. 10 (1976), n. 3, 70--71 (in Russian);
Berezin F., Laplace--Cazimir operators on Lie supergroups.  The
general theory.  Preprints ITEPh 77, Moscow, ITEPh, 1977; Berezin F.
{\it Analysis with anticommuting variables}, Kluwer, 1987
\bibitem[Bu]{Bu}
Bourbaki N., {\it Alg\'ebre commutatif}.  Ch.  V--VII,
Masson, Paris, 1985
\bibitem[L]{L}
Leites D., ed.  {\em Seminar on supermanifolds}, Reports of Stockholm
University, 1989-91, vv.1--34, 2100 pp.  (An expanded version of the
book: Leites, {\em Supermanifold theory}, Petrozavodsk, Karelia Branch
of the USSR Acad.  Sci., 1983, in Russian)
\bibitem[M]{M}
Macdonald I., {\em Symmetric functions and Hall polynomials}, Oxford
Univ. Press, 1995
\bibitem[S1]{S1}
Sergeev A., Laplace operators and representations of Lie
superalgebras.  In: [L], v.  23
\bibitem[S2]{S2}
Sergeev A., The tensor algebra of the identity representation as
a module over Lie superalgebras $\fgl (p, q)$ and $\fq (n)$. Mat. Sb.,
1984, 123, 3, 422--430
\bibitem[S3]{S3}
Sergeev A., Irreducible representations of solvable Lie superalgebras. In:
[L], v. 23
\bibitem[Wi]{Wi}
Weil A., Th\' eorie des points proches sur les vari\' et\' es diff\'
erentiables.  (French) G\' eom\' etrie diff\' erentielle.  Colloques
Internationaux du Centre National de la Recherche Scientifique,
Strasbourg, 1953, pp.  111--117.  Centre National de la Recherche
Scientifique, Paris, 1953
\bibitem[Wy]{Wy}
Weyl H., {\em Classical groups.  Their invariants and
representations}.  Princeton Univ.  Press, 1947
\end{thebibliography}
\end{document}